\documentclass{tMAM2e}

\usepackage{algorithm}
\usepackage[noend]{algpseudocode}
\usepackage{epstopdf}
\usepackage[bookmarksnumbered,colorlinks,bookmarks,citecolor=blue,linkcolor=blue,breaklinks,linktocpage]{hyperref}
\usepackage[all]{hypcap}%
\usepackage{subfigure}

\usepackage{diagrams}
\usepackage[all]{xy}
\usepackage{tikz-cd}
\usetikzlibrary{arrows}
\tikzset{commutative diagrams/.cd,arrow style=tikz,diagrams={>=latex'}}
\usepackage{graphicx}
\usepackage{multicol}

\usepackage{soul}
\usepackage{todonotes}

\theoremstyle{plain}
\newtheorem{theorem}{Theorem}[section]

\newtheorem{conjecture}[theorem]{Conjecture}

\theoremstyle{definition}

\theoremstyle{remark}

\begin{document}



\title{Quantum GestART: Identifying and Applying Correlations between Mathematics, Art, and Perceptual Organization}

\author{
\name{Maria Mannone\textsuperscript{a}$^{\ast}$\thanks{$^\ast$mariacaterina.mannone@unipa.it},
Federico Favali\textsuperscript{b},
Balandino Di Donato\textsuperscript{c},
Luca Turchet\textsuperscript{d}}
\affil{\textsuperscript{a}Department of Mathematics and Informatics, University of Palermo, Italy;
\textsuperscript{b}Universidad Nacional de Tres de Febrero, Buenos Aires, Argentina; \textsuperscript{c}Informatics Department,
University of Leicester, UK; \textsuperscript{d}Department of Information Engineering and Computer Science, University of Trento, Italy}
}

\maketitle


\begin{abstract}
Mathematics can help analyze the arts and inspire new artwork. Mathematics can also help make transformations from one artistic medium to another, considering exceptions and choices, as well as artists' individual and unique contributions. We propose a method based on diagrammatic thinking and quantum formalism. We exploit decompositions of complex forms into a set of simple shapes, discretization of complex images, and Dirac notation, imagining a world of ``prototypes'' that can be connected to obtain a fine or coarse-graining approximation of a given visual image. Visual prototypes are exchanged with auditory ones, and the information (position, size) characterizing visual prototypes is connected with the information (onset, duration, loudness, pitch range) characterizing auditory prototypes.
The topic is contextualized within a philosophical debate (discreteness and comparison of apparently unrelated objects), it develops through mathematical formalism, and it leads to programming, to spark interdisciplinary thinking and ignite creativity within STEAM.
\end{abstract}



\begin{keywords}
Gestural similarity, Gestalt, diagrams, Dirac notation, sonification
\end{keywords}
\begin{classcode}\textit{2010 Mathematics Subject Classification}:  81-00; 18C10; 42-00  \\
\end{classcode}

\section{Introduction: philosophical background}\label{introduction}

Understanding things that surround us in daily life also includes the ability to classify and even ``simplify'' them. This means collecting different things (objects, images, shapes of living beings) within classes or categories. This also means  recovering similarities between objects and between their transformational processes. Mathematics, and especially some branches of abstract mathematics such as category theory, provide the needed formalism to compare objects via transformations, and transformations between transformations.
Among the things around us, there are also sounds, and, as a subset of sounds, music. Musical structures also contain objects --- notes and recognizable themes --- that can be transformed.

Is it possible to find a key, a technique, to connect not only objects that easily appear as similar but also objects whose similarities are hidden? Is it possible to simplify an object into essential lines and to translate these lines from apparently unrelated objects, such as a tree and a symphony, a book of geometry and a painting? Can we imagine a world of ``visual prototypes'' and a world of ``sound prototypes,'' that can be eventually related within organized structures? Unity within the variety of human thinking can remind us of Leonardo da Vinci and of the Italian Renaissance.
While researching shared structures between sounds and images, and sounds and shapes/forms from nature, would it be possible and useful to apply some principles of a cognitive theory that specifically investigates forms and their organization?
To this aim, we can consider ideas from the Gestalt, a theory of the form seen as a whole that is more than the sum of its individual parts. This applies to visuals, but also to sounds \citep{gestalt_music, lerdhal}. Gestalt is related to perceptual organization. We may wonder if, in the arts, the generic ``form'' can be an ``artistic idea'' that is embodied in this or that  specific artistic code.
For this reason, we can ask: how are some cognitive topics, such as Gestalt, relevant to a mathematical theory of the arts? 

The search for hidden connections is reminiscent of contemporary philosophical approaches.
It is well known that the French philosopher Gilles Deleuze\footnote{Deleuze's contribution also involves diagrammatic thinking. According to \cite[page 43]{deleuze_}:  ``Le diagramme ne fonctionne jamais pour repr\'{e}senter un monde pr\'{e}existant, il produit un nouveau type de r\'{e}alit\'{e}, un nouveau mod\`{e}le de v\'{e}rit\'{e},'' that is, ``A diagram doesn't work to represent a world as a pre-existing one; it rather produces a new type of reality, a new model of truth.''} (1925-1995) developed the concept of {\em rhizomatic thought} \citep{hillier}. According to Deleuze, as the rhizome has many branches, it is possible to link different topics from different fields and study the ``hidden'' relations among them. Thus, several concepts can refer to a common concept. In this way, we can link together several seemingly different fields.
This is coherent with the word, ``quantum,'' from our proposed title which also refers to the discretization of rizhoma; ``Gestalt'' highlights the need for an organization and an identification of groups of rhizomes, and ``ART'' highlights the importance of this topic for a new approach towards arts. For this reason, we can coin a new expression, {\em Quantum GestART.}\footnote{The term ``GestART'' has been used by M. Mannone in a European application (H2020) in Summer/Fall 2017, and then used in this framework. The same word has been independently used in a 2018 online project: \url{https://www.univeur.org/cuebc/index.php/it/gestart-exhibition\#1}.}

Rhizomatic systems are mainly horizontal, as opposed to the verticality of Hegelian categorization. Rhizomas can be considered not only as static objects, but as events in space and time, that can change, can be composed together, and grow. This could give an idea of a comprehensive Gestalt for these topics. Also, creativity lies in the combination of rhizomatic entities, with the creation of new connections and new arrows. In this way, new structures and new hierarchies between elements can appear.

An overall form can be seen as the result of an evolution in space and time --- and this is the perspective considered in ``On Growth and Form'' by D'Arcy Wentworth Thompson \citep{darcy, lindenmayer, lindenmayer2}:
\begin{quote}
Organic form itself is found, mathematically speaking, to be a function of time... We might call the form of an organism an {\em event in space-time}, and not merely a {\em configuration in space}.
\end{quote}
Also, according to \cite[page VIII]{halle}, quoted in \citep[page 29]{lindenmayer2},
``The idea of the form implicitly contains also the history of such a form.'' Considering a visual form as (the result of) a collection of events in space and time may lead to its sonification, i.e., translation into sounds \citep{hermann2011sonification}. In fact, such a sequence of visual events can be mapped into a sequence of musical events in time. If we consider a little branch \citep{mannone_dissertation}, we can build it little by little in time as a temporal representation of the growth process. 
Thus, the sonified little branch would be a sequence of events in the auditory domain rather than the visual, which then aims to convey to the listener the idea of an organic shape that is replicated and grows.  One example is a short theme which is repeated and transformed. 

Reading an image as a collection of events implies its discretization. Discretization via Deleuzian rhizomas can be metaphorically compared with discreteness proposed in quantum mechanics. Thus, an original image can be thought of as a continuous form, and its decomposition onto a set of simple shapes. This idea will be described in detail in Section \ref{mathematical details}. The image discretization is provided by a superposition of rhizomatic entities. Gestalt is given by the collection of these entities --- it may be metaphorically compared with their categorical colimit\footnote{{\em Limits} and {\em Colimits} are concepts from category theory \citep{maclane}. As explained by \cite[page 12]{mannone}: ``Intuitively, {\em limits} and {\em colimits} are generalizations of {\em products} and {\em coproducts}, respectively. The product is a special case of the limit, with discrete categories. The coproduct (also called {\em sum}) is the dual of the product, obtained by reversing the arrows. Given an object $P$ and two maps $p_1:P\rightarrow B_1,\,p_2:P\rightarrow B_2$, $P$ is a product of $B_1,\,B_2$ if for each object $X$ and for each pair of arrows $f_1,\,f_2$ we have only and only one arrow $f:X\rightarrow P$ such that $f_1=p_1f,\,f_2=p_2f$ \citep{lawvere}, see diagram \eqref{product}.''
\begin{equation}\label{product}
\begin{footnotesize}
\begin{diagram}
&	&&		\mbox{X }		& &\\
&	&\ldTo(2,4)^{f_2}&\dTo	 ^{f}&			\rdTo(2,4)^{f_1}	& &\\
&	&&			P	& &\\
&			& \ldTo^{p_2}&&\rdTo^{p_1}		 &  \\
&			B_2&			&	& &B_1 &    \\
\end{diagram}
\end{footnotesize}
\end{equation}} \citep{maclane}. Why are we considering  the concept of colimit in such a framework, even metaphorically? Intuitively, let us imagine to have several arrows starting from each one of these rhizomatic entities, and reaching a form. A form that is well-recognizable, made up by rhizomatic entities as a particular building or a clear form made up of bricks. We may think of several LEGO$^{TM}$ bricks bound up to recreate a specific, easy-to-recognize, form. And this construction is ``universal,'' in the sense that the arrows starting from other bricks can be composed with the other arrows to reach the same form, and that composition of arrows is unique.\footnote{This ``unicity'' is a rough example of the ``universal property'' shown by some constructions in category theory, such as colimits and limits, dual of the colimits \citep{maclane}.}

The step of the image discretization via a ``base'' decomposition is represented by the step $F$, filter, in the mathematical procedure we propose in Section \ref{mathematical details}.
Then, we sonify the decomposition (step $T$). Finally, we can refine the obtained sequences to compose a complete musical piece (step $R$), via an interpolation of the obtained musical sequences. In some sense, from the discrete we get back to the continuous.\footnote{The idea that a shape may contain its story can be also compared with evolution, creating an ideal connection between the evolution of natural forms and the evolution of musical forms. This fact may be contextualized in a wider framework, within the evolution of Western (artistic) thinking.}

In this article, we summarize some ideas of recent literature about these questions. Then, we propose a theoretical development: a complex image can be filtered and approximated via the superposition of a set of schematic shapes with opportunely chosen coefficients. 
A minimal number of simple shapes -- that we will call ``visual vectors,'' or ``visual kets'' in Section \ref{method__} --  is needed to make the image recognizable, as in Gestalt principles. In fact, we use a formalism inspired from the bra--ket Dirac notation used in quantum mechanics \citep{dirac_ket}. 

An image can be reduced to a collection of vectors with coefficients. Then, the ``base'' can be changed from visuals to sounds: each considered set of visual vectors is associated with a set of sound vectors, a set that is  opportunely chosen according to former studies about music and sound. The new superposition of basic sounds and sound sequences gives a larger musical sequence. With an abuse of notation that we will explain later in this article, we can talk of a visual ``basis'' (or base) and of a sound ``basis'' (or base).
By changing the basis, we make a transformation from the visual domain to the sound domain.
Because a Fourier transform may be seen as a change of basis, we can metaphorically describe this process as a Fourier-like transform (and thus we will be using a $T$ symbol) from visuals to sounds. One of the goals of this study is to compose music. For this reason, the final result of the described process can be refined by a (human) composer and used as a complete musical piece.

This strategy can also be contextualized within the framework of sonification studies \citep{hermann2011sonification}. The resulting music will have its own Gestalt, that would be the sound rendition of the original visual Gestalt \citep{mannone_aesthetics}. If the final result can reasonably be brought back to the original image --- or at least to an equivalence class of images that contains the original one --- we can talk about some degree of {\em conservation} of the overall Gestalt, in a way that is inspired by conservational and invariance principles in physics.
Due to the generality of this strategy, the initial image can be a form from nature, an abstract collection of points and lines on a painting, and so on.

Our study is also deeply related to {\em crossmodal correspondences}. 
This phenomenon refers to people's tendency of matching information presented in one sensory modality to the information presented in another one \citep{spence}. 
Crossmodal correspondences are statistically--based experiences, and they are considered different from {\em synesthesia},\footnote{The {\em synesthesia} concept is familiar from poetry: let us recall the ``black scream'' in ``[...] All'urlo nero della madre / che andava incontro al figlio / crocifisso sul palo del telegrafo'' from the poem {\em Alle fronde dei salici}, in {\em Giorno dopo giorno} (1947) by Salvatore Quasimodo, Literature Nobel prize in 1959.} where the association is involuntary and more individual-dependent \citep{deroy_spence}. In some cases, the association between essential properties of sounds and essential properties of visuals allows the identification of the so-called {\em audio-visual objects} \citep{kubovy}; this seems to be relevant also in an artistic framework, where we may characterize multisensory artworks. In this paper, we will be dealing with the perceptual organization of different sensorial stimuli, that can be thought of as a visual Gestalt and an auditory Gestalt. We will refer to joint studies between Gestalt and crossmodal correspondences \citep{spence_Gest}. 
The analysis of visual and auditory Gestalt can also involve sounds in speaking and visual representations (see \citealt{nobile} for Italian language).    

The structure of the paper is the following: in Section \ref{mathematics and cognition}, we contextualize our work within the current interdisciplinary literature. In Section \ref{mathematical details}, we present the mathematical formalism. In Section \ref{pseudocode}, we show some sketches of a pseudocode, discussing how it would be possible to implement these ideas and strategies via a machine learning application.
In Section \ref{experiment}, we discuss future developments of the research, included the possibility of a cognitive experiment that would give concreteness to the theoretical approach we present here.
In the Conclusions, we consider additional philosophical and computational developments.
The method we describe in Section \ref{mathematical details} allows for translations image-to-music; however, reversing the arrows, music-to-image is also possible. This gives a precise and formal background to former approaches to music-image translations. Images include tridimensional objects whose space coordinates are mapped into time, loudness, and pitch \citep{mannone_libro}.

\section{Background: Cognition, Physics, and Mathematics}\label{mathematics and cognition}


\subsection{Cognition: Gestalt and crossmodal correspondences}


In this section we focus on examples of crossmodal correspondences between sounds and images, on Gestalt in music, and on Gestalt joined with crossmodal correspondences. According to \cite{spence} and \cite{spence_3}:
\begin{quote}
Crossmodal correspondences are defined as tendencies for a certain sensory feature (or dimension) to be associated or matched with another feature (or dimension) in a distinct sensory modality.
\end{quote}
Several examples of them have been found between vision and audition. For instance, we can cite the widely shared association between an ascending melodic line and an ascending vertical movement \citep{spence, zbi}. Literature describes crossmodal correspondences between visual, auditory, tactile, and even olfactory stimuli  \citep{spence_3}. We can wonder about the existence of general or generic sensorial concepts that can be embodied in different, specific senses. 
Crossmodal correspondences, according to \citep[page 885]{spence_3}, ``is something that occurs between {\em concepts}.'' If we extend this idea to organized stimuli, the question about the ``translatability'' of an artwork from a sensory field to another would be meaningful. In our framework, translations are made via category-inspired diagrams. In fact, the mathematics we use to investigate these problems, category theory, has a high abstraction power. In particular, the search for abstraction and  invariants (see the end of Section \ref{mathematical details})
is also outlined in cognitive literature such as the work of  \cite[page 62]{gestalt_music} which states, ``The process of abstraction [...] is characterized by the extraction of invariants.''

The second topic of cognitive sciences relevant to our study is Gestalt theory. According to \cite{gestalt_music}, besides its intuitiveness for visuals, Gestalt laws often need an experimental validation to be extended in other areas. Reybrouck also points out ``a fundamental isomorphism (in a metaphorical sense) between the structure of the object and the psychophysiological structure  (Gestalt-like patterns of excitation in the brain) of the perceiver.'' This highlights a correspondence between what we can call the ``main lines'' characterizing the structure of an object, and the ``patterns'' in the brain that are activated when these main lines match with them. As we will see in Section \ref{mathematical details}, this idea\footnote{Curiously, \cite{gestalt_music} explicitly refers to algebraic methodologies and cites the {\em functor} term, but it appears as more in the sense of a function than of a categorical functor. A functor is a morphism between categories \citep{maclane}.} can support a decomposition of visual things into a superposition of base elements, and Gestalt
as their metaphorical colimit.\footnote{Gestalt as a colimit: all simple elements are connected into a colimit. We are considering here colimits and not limits because arrows are directed from small entities to the complete, recognizable form, the Gestalt. As explained above in our example with LEGO$^{TM}$ bricks, let us consider a collection of simple shapes that can be combined together to obtain a recognizable form. Thus, the arrows we are considering are directed from the simple shapes to the `complete' form. A reviewer thought of a colimit of `visual things' as patterns in the brain. This is a very interesting remark that can be addressed in future research; also, there are recent studies between categories, brain, and consciousness, see \citep{hayato}.}



Because music happens in time, a Gestalt approach to music should be considered as a dynamic process,
and music should involve a ``description in terms of processes''  
\citep{gestalt_music} rather than a description of static structure.\footnote{In the approach proposed in \citep{mannone_dissertation}, images are seen as the result of processes, thus the use of Gestalt from images to music in this framework appears as more coherent.} According to \cite[page 64]{gestalt_music}, ``In contrast with a geometrical figure, that is described as a whole when looking at it, a musical figure needs a successive presentation.'' However, if we consider a geometrical shape not -- yet -- as a whole, but as a ``making'' of it in terms of movements, we can easily connect it with a sequence of musical gestures (see Paragraph \ref{gest_generator_paragraph}).
If we can see images in their making, we can also do the opposite: we can try to discretize music in terms of short sequences that are perceived as a whole.
According to \citealt{gestalt_music} (page 64):
\begin{quote}
The grasping of its meaning [music], then, is polythetic in that it can be grasped only by hearing the music as it unfolds step by step (Wright, 1995). This hampers the holding of a musical Gestalt as an immediate and directly experienced whole. Music, however, can be grasped in a monothetic way, if the discursive processes are coded as discrete things as in the case of conceptualization (McAdams, 1989).
\end{quote}
This seems coherent with the quantum approach proposed in Section \ref{mathematical details}, where simple visual units are connected with simple musical units, and groups of visual units are thus mapped into groups of musical units.
Groups of visual units are meaningful entities that can be considered as a {\em visual Gestalt}. Equivalently, groups of musical units can be considered as a {\em sound Gestalt}. Thus, mapping groups of visual units into groups of musical units can ideally be seen as a mapping from a visual Gestalt to a sound Gestalt. We wonder if such a mapping can be possible and effective.

\cite{spence_Gest} summarizes recent research on Gestalt and crossmodal correspondences. He also considers open questions such as ambiguities of perception in Gestalt and their possible crossmodal equivalent. Another question about the possible influence of crossmodal perceptual grouping on intramodal grouping, that is the opposite of the first question. Experimental evidence shows that there is an influence of a grouping within a sensory area towards a crossmodal sensory grouping which involves that area and another one \citep{spence_Gest}. The first intersensory Gestalt that has been validated is the audiotactile musical meter perception.\footnote{This can easily be contextualized within studies about music and gestures in terms of physical movements. If we interpret visuals in terms of gestures, our quest for Gestalt translation from visual domain to auditory domain via mathematics and gestures may appear as well-defined. Also, according to \citealt[page 11]{spence_Gest}, crossmodal correspondences ``can be conceptualized in terms of cross-modal grouping by similarity.''}

Thinking of categories, we can see crossmodal correspondences as arrows between different sensory areas --- e.g., visual and auditory areas --- and Gestalt as arrows that start from a collection of elements on each sensory area and reach a point (as a  categorical colimit), the abstract form. We wonder if we can draw arrows, that is, if we can connect a Gestalt of a sensory area (e.g., Gestalt of image) with a Gestalt of another sensory area (e.g., Gestalt of sound), represented, respectively, as two ``final'' points in diagram \eqref{diagram_0}, where the question mark indicates a connection not verified yet.
\begin{equation}\label{diagram_0}
\begin{diagram}
\mbox{Gestalt (image)} & \rTo^{?} & \mbox{Gestalt (sound)}
\\ \uTo^{extraction} \uTo & & \uTo^{extraction} \uTo
\\ \mbox{image} & \rTo^{crossmodal} & \mbox{sound} 
\end{diagram}
\end{equation}
The reference to Gestalt as an organization principle and crossmodal correspondences as comparisons between different sensorial stimuli make us think of a superior, more general, ``suprasensory'' organizing principles \citep{spence_Gest, allen_kolers}, or of trans-modal Gestalts \citep{yu}.
In particular, \cite{yu} consider additivity of groups that is relevant for models with states superposition.

These ideas can be connected with a recent development of the mathematical theory of musical gestures: the concept of {\em gestural generator} \citep{mannone}, see Paragraph \ref{gest_generator_paragraph}, that can be in some way compared with these approaches. Intuitively, a simple drawing and a short musical sequence can be seen as {\em similar} if they appear as being produced (generated) by the same movement, with pencil on paper, or with hands/breath/bow on musical instruments \citep{mannone}; see diagram \eqref{gest_1}.
\begin{equation}\label{gest_1}
\begin{diagram}
  & & \mbox{gesture}  & &
\\ & \ldTo^{generates} & & \rdTo^{generates}
\\ \mbox{image} & & \rTo^{similarity} & & \mbox{music}
\end{diagram}
\end{equation}
 This conjecture has been experimentally validated for simple shapes and simple musical sequences \citep{collins, mannone_collins}. Our theoretical framework could be used to contextualize other experiments about sounds and gestures, even with non-musical sounds or sounds not derived from sonification but just noises connected with the drawing production. In particular, we can refer to a recent work on synchronization of gestures with synthetic friction sounds \citep{thoret}. Diagrammatic thinking can help investigate mutual influences and biases between audition, gestures, and drawing production. This is the case of sounds variations that may induce people to draw ellipses instead of circles \citep{thoret3}, and the perception of biological movements via sound detection \citep{thoret2}.\footnote{In these works, the $\frac{1}{3}$ law is cited, that is, the connection between timbre variations and velocity profiles. This reminds one of the $\frac{2}{3}$ power law describing the connection between the kinematics of handwriting and the trajectory of a movement, between the curvature of handwriting and the pentip angular velocity \citep{lacquaniti, plamondon}.}
The hypothesis of a connection between different senses and media towards the representation of ``higher'' and sense-independent concepts, may be supported by recent research in the so-called ``supramodal brain'' \citep{rosenblum}. According to \citealt[page 65]{rosenblum}, within the supramodal conception of the perceptual brain, ``task, rather than sensory modality, is the primary organising principle.'' Also, research focused on speech shows that ``supramodal speech mechanism also seems sensitive to haptic and kinaesthetic input'' \citep[page 73]{rosenblum}. This may be related  with the idea of gestural generator; see Paragraph \ref{gest_generator_paragraph}.
 
Let us also consider the possibility of translations, transfer of Gestalt from one domain to the other. About crossmodal Gestalt, we may wonder about the existence of some {\em uncertainty principle}; we may not be able to discern with the same precision all groupings and so on in one medium, e.g., in visuals, and, at the same time, to perfectly recognize groupings on another medium, e.g., sound, if performed simultaneously. Formally, we have the inequality of Eq. \eqref{uncertainty}.
\begin{equation}\label{uncertainty}
\begin{split}
&\mbox{(incertitude in finding groups/similarities in visuals)} \\&\times\mbox{(incertitude in finding groups/similarities in music)}
\\ &\geq\mbox{(constant)}
\end{split}
\end{equation}
This means that either we hear with extreme attention or we make a visual analysis with extreme attention. An intermediate value of attention should be supported by artistic realizations such as {\em Optical Poem} (1938) by Oskar Fischinger. Already \cite{gabor} proposed an uncertainty principle for sound (only), with a formalism borrowed from physics. The inequality of Fig. \eqref{uncertainty} is not directly inspired of non-commutativity, but of precision levels. However, we can draw a diagram of attention levels starting from diagram \eqref{diagram_0}, and we can experimentally investigate its commutativity. The inequality relationship of Eq. \eqref{uncertainty} implies the use of short temporal windows, equivalent to the sequences of Fig. \ref{basis} for each small group.
In conclusion, the relations between cognition and neural system can be presented more formally by using also a mathematical approach to the neuro-cognitive system, based on ``dynamic'' category theory \citep{ehresmann_van, kainen}. In this model, called ``Memory Evolutive Systems,'' developed by Andr\'{e}e Ehresmann and Paul Vanbremeersch, the mental objects are colimits of diagrams that represent synchronized neural groupings. According to \citealt[page 228]{kainen}: ``Its creators, Ehresmann and Vanbremeersch, have proposed a unique synthesis of various ideas and concepts, organized around the notion that life is an evolving hierarchical category with a functorial dynamic.'' Thus, this idea also connects cognition with natural studies.
This approach led to very recent developments about mental representations investigated via categories \citep{hayato}.

\subsection{Physics: quantum mechanics}

References to cognition may be relevant for a mathematical theory of the arts for two reasons. First, as it happens in physics, mathematics can give more solutions than the ones that are really coherent with natural phenomena. Experiments are needed to verify the consistency of a theory with the physical reality. Likewise, cognitive experiments can be useful to validate specific approaches or conjectures of mathematical theories of the arts. Second, cognitive-relevant topics may give new blood to theoretical advancements, 
that stimulate the development of new theories or the generalization of existing ones.

The parallel between physics and mathematical approaches to the arts is not only a metaphor. To investigate cognition, visuals, and sounds, there have been attempts to exploit ideas and methods from quantum mechanics, using the main concepts of quantization, destructive measure, and probability \citep{gabor, khrennikov, rocchesso_mannone}. 
The step from continuity to quantization can sometimes make evaluations easier, and it can allow analysis in a new light. A renowned example is about sound itself. Acoustical waves are often used as a paradigm for quantum waves \citep{gabor} 
conversely, if we use the formalism of quantum waves to describe sound, we can also borrow the duality wave-particle and introduce quanta of sound \citep{gabor}. Other examples include analysis of images \citep{youssry}, psychology \citep{khrennikov}, musical forces and cognition \citep{graben, graben2}, quantum information, cognition, and music \citep{dallaChiara}, adaptation to music of techniques to evaluate memory in quantum systems \citep{mannone_compagno}, and vocal theory based on quantum formalism \citep{rocchesso_mannone}, where vocal categories are used as probes to describe sounds in general.

\subsection{Mathematics: category theory and diagrammatic thinking}

Another topic often applied to music is a mathematical theory: category theory. A category is constituted by objects and morphisms between them, verifying associative and identity properties \citep{maclane}. This is a general theory that starts with the very basic concepts of objects (visually represented as points), transformations between them (visually represented as arrows), and transformations between transformations (arrows between arrows). We can see an entire category as a point, and we can draw arrows between categories (called {\em functors}), creating nested structures.
Finally, a mathematical tool we will use is a Fourier-inspired transform, which will be defined in Section \ref{mathematical details}. The use of Fourier transforms not only for sound synthesis but also for music analysis has been widely exploited in music theory, to characterize rhythmic or harmonic distributions \citep{amiot_book}.

Categories play an important role in attempts to unify different branches of mathematics itself \citep{caramello}. Also, they have been used to investigate human cognition \citep{phillips}, neuroscience \citep{ehresmann}, aesthetics, art production, and art evaluation processes \citep{kubota}\footnote{Curiously, this paper has been submitted on the same day of a paper on categories and composer/conductor/listener, see \citep{mannone}.}, the aesthetics of mathematics \citep{mannone_aesthetics}, and physics \citep{baez}. 
Categories have also been applied to music \citep{mazzola, popoff, jed, jed2, fiore, mannone}. We can describe, with categories, musical structures as elements and transformations between them, musical gestures \citep{mazzola_andreatta,arias}, recursive musical gestures \citep{knots}, and even musical instruments \citep{mannone_favali}, and music and dance \citep{mannone_turchet}.
In recent literature, there have been some pioneering connections between quantum mechanics and diagrammatic thinking \citep{coecke}.

Category theory is the general framework of recent studies we are starting from. However, in this paper we will not make use of precise categories; we will rather consider a more generic diagrammatic thinking as a visual aid to describe transformations and transformations between transformations.

\subsubsection{Category theory and nature}

The categorical concept of ``transformation between transformations'' can be applied to shapes from visual arts, and even to shapes and forms from nature. For example, we can compare the growth of a fish with the growth of another fish via an arrow between arrows \citep{mannone_libro_nuovo}. Categories in biology had already been applied to nature in terms of input-output processes \citep{letelier}. Other relevant mathematical approaches to shapes from nature are the classic studies by D'Arcy Wentworth Thompson \citep{darcy}, that investigates the geometry of nature, and the studies of Aristid Lindenmayer \citep{lindenmayer}, that investigates the growth of plants. These studies are relevant in particular for computer graphic applications. D'Arcy W. Thompson and Lindenmayer's studies do not involve category theory, but they can be re-read in its light \citep{mannone_libro_nuovo}.
According to \cite{gestalt_music}, references to {\em processes} and {\em growth}, jointly with the concept of {\em form}, suggest an immediate connection between Gestalt and forms and growth of natural forms \citep{darcy}.
The categorical formalism has been applied to processes in biology  (not to morphology) by Robert \cite{rosen}, and extended to the emergence of cognition from neurons by \cite{ehresmann_van}.
Finally, classification can also be investigated in light of categories, by using the concept of colimit to connect species, genera, and families with composition of arrows leading to the same ``essential idea'' \citep{mannone_libro_nuovo}.

\subsubsection{Gestural generator}\label{gest_generator_paragraph}
We can intuitively define a ``gestural generator'' via its action as we will see.
The concept can be applied within music, and between music and visuals. Within music, let us consider two performance gestures: a movement to play a note with an unspecified dynamics (loudness level) on a piano, and the movement to play a note with an unspecified dynamics on a vibraphone, for example. To play a note with dynamics {\em forte}, we have to ``deform'' these movements. The two {\em forte} movements, on the piano and on the vibraphone, will then be {\em similar} with respect to the loudness. {\em Similar} changes in piano and violin loudness (and timbre) are provoked by a {\em similar} {\em gestural generator}, that operates a similar change in the respective gestures --- e.g., acceleration increase for hands/hammers, or of airstream pressure for a wind instrument.\footnote{Diagram \eqref{generator}, from \citep{mannone}, shows a mapping from an abstract collection of points and arrows ($\Delta$) to the topological space $\vec X$ (containing curves in the $X$, the space of points in physical space and time). In the diagram, $g_{\phi}$ represent the physical ($\phi$) gesture without any specified dynamics, while $g_{\phi}^F$ is the gesture with the dynamics {\em forte} (F). The double arrow $F$ represent the deformation of the unspecified gesture into the ``forte'' gesture. This is a specific {\em forte} operator. We can think of the {\em ``forte'' gestural generator} as a generic ``forte'' gesture, to be adapted in the specific space of gestures for each musical instrument that deforms anonymous movements into movements that produce {\em forte} loudness on musical instruments.} Of course, all the changes must be made in the appropriate spaces. Striking a flute with a hammer will not produce a {\em forte} flute sound: instead, it will just break the flute! A formal definition of gestural similarity is given by \cite{mannone}. We report the Heuristic Conjecture about gestural similarity of musical gestures \citep[page 22]{mannone}:
\begin{conjecture}
Two gestures, based on the same skeleton,\footnote{A ``skeleton'' is an abstract structure of points and arrows; see $\Delta$ in diagram \eqref{generator} as an example. Such a generic structure is then specified into a collection of points and arrows in a topological space. In our case, we consider physical space and time.} are similar if and only if they can be connected via a transformation:
\begin{enumerate}
\item that homotopically transforms one gesture into the other, and
\item that also leads to similar changes in their respective acoustical spectra.
\end{enumerate}
Homotopy is a necessary, but not sufficient, condition to get similar gestures.
\end{conjecture}

While dealing with music and images, we can think of a gesture that can either generate an image or a drawing: a staccato sequence is similar to a collection of points, or a legato phrase to a continuous line, and so on. Thus, pairs of images and sounds can share some degree of {\em gestural similarity} if their appear as being generated by the same movement. See diagram \eqref{gest_1} and its development in Section \ref{mathematical details}, with diagram \eqref{gest_diagram}. This approach can help investigate interactions between gestures, images --- and in particular, drawing --- and sounds --- and in particular, music. We can re-read in this light the pioneering work by Alexandre Truslit \citep{repp} including the comparison of landscapes and melodic profiles of folk songs, and the focus on ``inner motion'' and movement in music analysis. This analysis can be extended to contextualize experimental studies about gesture production as responses to auditory stimuli and in particular music \citep{kussner}, and the influence of vision of gestures in music perception \citep{lipscomb}.
\begin{equation}\label{generator}
\begin{tikzcd}
\Delta
  \arrow[r,bend left=70, "g_{\phi}"{name=U, above}] 
  \arrow[r,bend right=50, "g^F_{\phi}"{name=D, below}]  
& 
\vec X
\arrow[Rightarrow, from=U, to=D, "F"]
\end{tikzcd}
\end{equation}


\subsection{The proposed method}\label{method__}

Mathematical formalism can thus be applied to music, images, and nature, and it can even help connect all these areas together. For example, we can first analyze a visual shape, derive a mathematical model, define a mapping strategy, and apply the strategy to make music.
Visual shapes do not intrinsically have a temporal dimension. However, we can imagine creating an image through a drawing gesture, and we can see the main lines and points of an image as a sequence of events (sounds) in time. A form that changes through time can be represented as a sequence of changing melodies, as the sonification of a series of pictures (as in chronophotography). 

In this paper, we propose to join hints from quantum formalism, diagrammatic thinking, cognition, and current research in mathematical music theory \citep{mannone}, as well as some insights from Gestalt theory applied to music and visual arts. 

Also, the Gestalt is useful when we look for the minimal decomposition: we can decompose a complex image in terms of simple shapes (points, lines, geometric figures), that we may call here ``basis vectors.'' However, this has a very specific meaning in algebra, being strictly connected with the concept of space dimensions (e.g., in a 2-dimensional space, a basis should contain two linear independent vectors).\footnote{The visual kets of a set can be considered as linearly independent if none of them can be obtained as a linear combination of the other ones; similarly for the sound kets. We can consider two kets $|\phi_i^v\rangle,\,|\phi_j^v\rangle$ as orthogonal if their (gestural) similarity degree is zero \citep{mannone}, and thus their scalar product is zero: $\langle\phi_i^v|\phi_j^v\rangle=0$. However, if they are identical, $\langle\phi_i^v|\phi_j^v\rangle=1$. Two kets with a similarity degree strictly between $0$ and $1$ have a scalar product between $0$ and $1$. Similarly for the sound kets. Thus, the scalar product assumes the meaning of a comparison between simple shapes or between simple musical sequences.} However, in the method proposed in this article, we do not need a number of basis vectors equal to the dimension  of the space: the number of vectors is instead connected with the precision of the approximation of the given form, and it depends on the complexity of the original image. Thus, we may need another name besides ``basis vectors.'' Because of the reference to Dirac notation, we can speak of ``visual kets'' while dealing with images, and of ``sound kets'' while dealing with sounds.
We would like to convey the idea of minimal decomposition of a form in terms of simple shapes, that is, in terms of some visual kets or of a few combinations of them. We need a minimal number of visual kets to make an image recognizable; and in the definition of some minimum principle, Gestalt theory intervenes.
Also, we can define a basis of sounds/sound sequences. Each sound sequence should correspond to a basic image via an intuitive sonification: we can imagine both visuals and sounds as being generated by the same gestural generator. This association is an application of the definition of gestural similarity \citep{mannone} and has connections with recent studies in musicology, about the connection between movement and sound \citep{zbi}.

The validation of the proposed development requires  cognitive experiments to establish what a better basis for images and sounds can be, 
as well as if, and to what extent, Gestalt can be preserved upon a basis change. ``Change of basis'' means the transition from visual to sound kets (and vice versa). Formally, we indicate the correspondence between each element of the considered set of visual kets (what can be called a ``visual basis'' with an abuse of notation) with each element of the considered set of sound kets (what can be imprecisely called a ``sound basis''), and each coefficient of the visual kets with each coefficient of the sound kets. In this way, the decomposition of an image (meant as a drawing or a painting, or even a sculpture) in terms of visual kets and visual coefficients, can be translated into a decomposition of sound kets and sound coefficients, constituting a rough musical piece. Such an approximative musical piece can be refined to obtain a real and complete musical composition. In fact, with ``basis change'' we indicate the transition from the domain of visuals to the domain of sounds. We were inspired by the fact that the Fourier transform from the time domain of a signal into the frequency domain can be seen as a basis change. In the following, we will be using the term ``transform'' just to indicate this metaphorical analogy.

In Section \ref{mathematical details}, we describe the mathematical formalism, which joins diagrams from category theory with bra-ket notation typically used in quantum mechanics.

\section{Bases, filters, and transforms}\label{mathematical details}

Let us give some mathematical details about the described procedure, starting with an intuitive description. The rough idea is the following:
$$ \mbox{musical work} = T(\mbox{visual work}),$$
where $T$ is a transform, which we will define soon.
With ``musical artwork,'' we can just indicate a musical piece, while with ``visual artwork'' we mean a drawing, a painting, or a sculpture, that can be translated into sound via $T$.

Formally, we can define a vector space for visuals, a vector space for sounds, and a transition between them. An image can be decomposed into a selection of simple visual shapes, and such a (visual) decomposition can be translated into a sound decomposition, living in the sound vector space. Such a sound decomposition constitutes a rough musical composition.

Let $V$ be the space of visuals (the space of visual shapes), and let $K$ be a field. Given $v_1,\,v_2, v\in V$ and $\lambda \in K$, we have: $V\times V\rightarrow V:(v_1,\,v_2)\mapsto v_1+v_2$, and $K\times V\rightarrow V:(\lambda,v)\mapsto \lambda V$. Similarly, let $S$ be the space of sounds (in the following examples, the space of musical sequences), and $K'$ another field. The visual vector space is given by the set of visual elements $V$, an operation of ``sum,'' $+$, that we suppose to be linear, and an operation of scalar product, $\cdot$. Let $u,\,v$ be two elements of $V$ and $\lambda$ a scalar of $K$; we have $+_V:(u,v)\rightarrow u + v \in V$, and $\cdot_V:(\lambda,u)\rightarrow \lambda\cdot u\in V$. Similarly, the sound vector space is given by a set of sound elements ---  here, short musical sequences --- $S$, an operation of sum $+_S$ and an operation of multiplication by scalars of $K'$ indicated with $\cdot_S$.
In this framework, the meaning of the addition is to superpose several simple shapes or several musical sequences within their respective spaces; the meaning of the multiplication by a scalar is to refine resizing and space positioning of visual or sound kets (the latest ones, in the space of time and sound parameters).

\begin{figure}
\centering
\includegraphics[width=13cm]{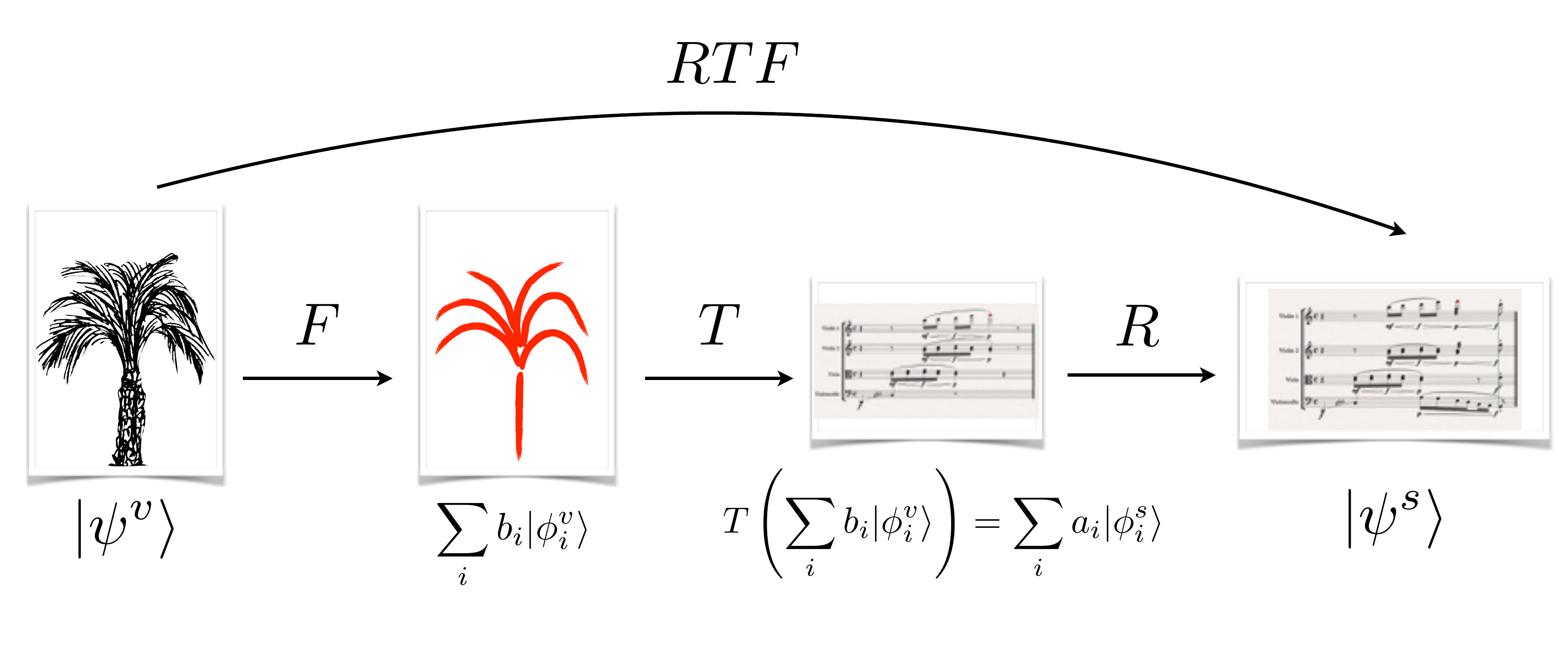}
\caption{The scheme RTF with the possible filtering (F), transform (T), and refinement (R) from a palm ({\em Butia capitata}) image. The diagram also shows the composition of the three operations. 
}
\label{RTF_image}
\end{figure}

In fact, the elements of the set of visuals are visual shapes, simple or complex. Let us consider original images ($|\psi^v\rangle\in V$) and visual kets as the elements living in the space $V$ (we can substitute $v_1,v_2,v\in V$ with $|\phi_1^v\rangle,|\phi_2^v\rangle,|\phi_3^v\rangle$, and we also have $|\psi^v\rangle\in V$), and the set of sound kets and complete musical pieces in $S$ (that is, $|\phi_1^s\rangle,|\phi_2^s\rangle,|\phi_3^s\rangle, |\psi^s\rangle\in S$, where the last one is a complete musical piece).\footnote{Formally, kets would require an Hilbert space, but here the imaginary component is zero or it is not considered. It could be identified with the mind component, that is embodied into the real part, that is, an artwork, during the artistic creation process \citep{mazzola}.}   The (visual) basis elements (that is, a selection of visual kets) are simple shapes, or, better, elementary shapes such as points, segments, and curved segments. Other simple shapes such as triangles and circles\footnote{In $\mathbb{R}^n$ a circle is not a vector but a set of such vectors.} can be obtained via a superposition of elementary shapes; e.g., a circle can be obtained with two essential lines: a pair of inverted curved segments, that is, two arches. For simplicity, when needed, in some points of the paper we will indicate as ``simple shapes'' also these simple superpositions of elementary shapes. 
Above, we already indicated the (essential) shapes with the visual kets $|\phi^v_i\rangle$. Similarly, we indicated with $|\phi^s_i\rangle$ the sound kets: elementary melodic profiles and articulations. In analogy with Euclidean vectors, visual (and sound) vectors have a magnitude (that we can assume unitary at first) and a direction (that is, in our research, the direction of time). In fact, musical sequences have a direction in time; simple visual shapes have a direction given by the direction of drawing; we can thus see visuals as events in space and time, analogously to what had been discussed about forms from nature \citep{darcy, halle}. The (initially unitary) magnitude and the position in the space of these vectors is modified by their coefficients. In the space of visuals, the coefficients --- indicated here with $b_i$ --- give information about size and position of the $i$-th basis elements. We can also think of ``rich'' coefficients telling us also how many kets of the $i$-th type are used. For example, in Figure \ref{RTF_image}, the image of the palm is simplified via the superposition of a segment and of several curved segments, corresponding to two kets: a segment and an arch. Coefficients $b_i$, thus, tell us how many kets are needed and where, in which points of the space, they have to be put to filter the original image in a recognizable way. Similarly, sound coefficients --- indicated here with $a_i$ --- tell us where (in time, thus, with which onset), for how long (with which duration), how loud, and how high in pitch the basic musical sequences have to be put. If we imagine a space of musical parameters with time, loudness, and pitch \citep{mannone_libro}, coefficients $a_i$ are positioning and resizing sound kets within that space. If we think of a generalized musical space with more parameters, we can dare to make more associations, for example between colors in $V$ and timbre in $S$, or contour details in $V$ and attack modes in $S$. For example, a vertical segment in the $V$ space may correspond to a cluster in the $S$ space, an horizontal segment in $V$ to a long, single note in $S$, and an arch in $V$ to a rising--lowering melody in $S$. How long the cluster or the melody has to last, and which precise pitches have to be played, and what loudness, is told by sound coefficients.
Even a single point has a meaning in $V$, and it can be translated in $S$ with a {\em staccato} single sound. The definition of ``good'' visual and sound kets can involve experiments, see Section \ref{examples}.

About the use of kets on a real field $K$ (while in quantum mechanics $K$ is complex), see \citep{sparavigna}.\footnote{\cite{sparavigna} consider kets as living in a subset of the vector space that can be of interest also for our study, both for the space of visuals and the space of sounds. It is required that: there is a zero vector $|0\rangle$; given $|a\rangle$ and $|b\rangle$ their sum $|a\rangle+|b\rangle$ is in the subset; the product by a scalar is also in the subset; for each $|a\rangle$ in the subset there is the opposite element $-|a\rangle$. In our study, the sum can be the superposition of shapes or musical sequences, and the zero vector the absence of a visual form (e.g., the empty visual space) and of musical sequences (e.g., the silence). The opposite $-|a\rangle$ can be the time-reverse of an image $|a\rangle $ --- that is, imagine starting with the finished simple drawing and going back to the empty space, or playing backward in time and going back to the silence. The operation of time reversing can be seen as some deleting operation, or, more interestingly, as just a time direction change: playing a melody from the latest to the first note, or drawing a simple shape from the last point to the first one.}

As a general remark, even if we are speaking of basis elements, there is no specific relation --- at least in this starting model --- with the dimension of the space. In fact, the more simple the original image $\psi^v$, the less visual kets $\phi_i^v$ are needed to unambiguously represent its Gestalt. Thus, the number of kets is connected with what we can call ``Gestalt-precision.'' See Figure \ref{ball} for an explanation.

Let us use the Greek letter $\psi^v$ to indicate the original image in $V$, and $\psi^s$ to indicate a complete, well-refined, musical piece.
The basis elements are however indicated with $\phi^v_i,\,\phi^s_i$, where $i$ indicates the $i$-th element in the considered decomposition. We are not using capital letters $\Psi$ or $\Phi$ here.

Let us start from the concept of superposition to develop a basic formalism, and from a complete image, e.g., a drawing or a picture.

First, the initial image (the ``visual work'') can be approximated via a superposition of (the least number of) simple shapes creating a sketch, such as a skeleton, and an overall envelope. Such a superposition makes the original shape recognizable, see eq. \ref{eq_1}, where $|\phi^v\rangle$ is the original, initial image, $b_i$ the coefficients for the image (visuals), and $|\phi_i^v\rangle$ the vectors for the visuals.  
\begin{equation}\label{eq_1}
|\psi^v\rangle \sim \sum_ib_i|\phi_i^v\rangle
\end{equation}

\begin{figure}
\centering
\includegraphics[width=8cm]{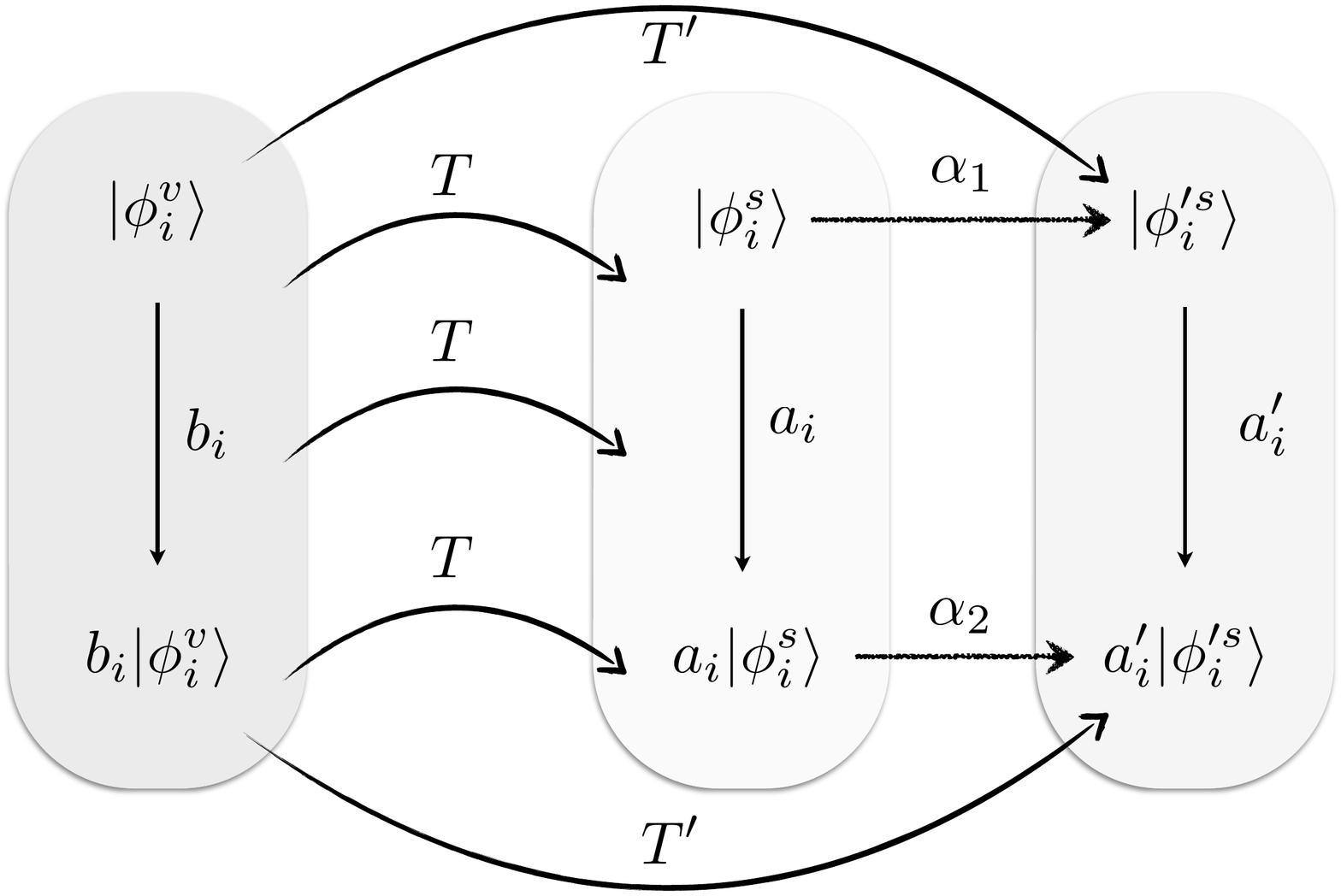}
\caption{Transform $T$ connects simplified images with sonified simplified images. A different sonification strategy, represented by $T'$, leads to different sonified images. The comparison between the action of $T$ and the action of $T'$ is indicated by transformations $\alpha_1$, $\alpha_2$.
The visual objects are the simple images indicated as kets $|\phi_i^v\rangle$, and their coefficients $b_i$ ``center'' each simple image in a point of the tridimensional space with a given size. Similarly, sound coefficients $a_i$ tell us in which point of sound space a sound ket $|\phi_i^s\rangle$ has to be put --- i.e., which time, loudness, duration shall we give it. Even if we are using diagrams and a category-inspired way to represent ideas, we are not precisely using a categorical formalism here: in fact, abstract kets $|\phi_i^v\rangle$, $|\phi_i^s\rangle$ are ontologically different from the same kets put in specific points of visual or sound space respectively, on a specific time with a given size or duration and loudness, that is, kets with coefficients: $b_i|\phi_i^v\rangle$, $a_i|\phi_i^s\rangle$.
}
\label{natural}
\end{figure}

Mapping from visuals to sounds is made via a basis change, see eq. \ref{eq_2} for coefficients and vectors ($a_i$ are the coefficients for the sound and $|\phi_i^s\rangle$ the vectors for the sound). We use index $v$ for visuals, and $s$ for sounds.
\begin{equation}\label{eq_2}
a_i|\phi_i^{sound}\rangle \iff b_i|\phi_i^{visual}\rangle
\end{equation}
Such a correspondence is used to evaluate transform $T$ in eq. \ref{procedure}. The correspondence between vectors for visuals and vectors for sound is established {\em a priori} (see gestural similarity for the ideas behind), while the value of coefficients and the number of basis terms involved depends case by case.  In Eq. \ref{procedure}, $|\psi^v\rangle$ indicates the original image, while $|\phi_i^s\rangle$ the sound kets, obtained via a mapping from the visual kets $|\phi_i^v\rangle$ used to filter $|\psi^v\rangle$. We might assume that $T$ is linear: $T(b_i|\phi_i^v\rangle+b_j|\phi_j^v\rangle)=T(b_i|\phi_i^v\rangle)+T(b_j|\phi_j^v\rangle)$.
However, $T(b_i|\phi_i^v\rangle)=a^i|\phi_i^s\rangle$, while linearity would require $T(b_i|\phi_i^v\rangle)=b_iT|\phi_i^v\rangle=b_i|\phi_i^s\rangle\neq a_i|\phi_i^s\rangle$. Thus, we can talk of {\em semilinearity}\footnote{Given two vector spaces $V, S$ over fields $K,K'$ respectively, $f:V\rightarrow S$ is a semilinear map if there exists a homomorphism $\sigma:K\rightarrow K'$ such that, for all $u,v\in V$ and $\lambda\in K$, we have $f(u,v)=f(u)+f(v)$ and $f(\lambda u) = \sigma(\lambda)f(u)$.} of $T$.
Also, we can suppose that $T$ can be represented by a diagonal matrix: $T_{ij}\neq 0\iff i=j$. In this way, we have $T(b_i|\phi_i^v\rangle)+T(b_j|\phi_j^v\rangle)=a_i|\phi_i^v\rangle+a_j|\phi_j^v\rangle$, and we can keep the correspondence between each visual ket with each sound ket. Of course, the condition of diagonality is a strong simplification, that can be useful for a first implementation.
\begin{equation}\label{procedure}
|\psi^s\rangle \sim \sum_ia_i|\phi_i^s\rangle = T\left(\sum b_i |\phi_i^v\rangle\right)\sim T(|\psi^v\rangle)
\end{equation}
Diagram \eqref{diagram_1} a commutative diagram contextualizes the whole procedure in a category-inspired framework.
\begin{equation}\label{diagram_1}
\begin{diagram}
\sum_ib_i|\phi_i^v\rangle & \rTo^{\sim} &  |\psi^v\rangle
\\ \dTo^T & & \dTo^T
\\ T\left(\sum_ib_i|\phi_i^v\rangle\right) & \rTo^{\sim} & T\left( |\psi^v\rangle \right)
\\ = & & =
\\ \sum_ia_i|\phi_i^s\rangle & \rTo^{\sim} & |\psi^s\rangle
\end{diagram}
\end{equation}
Let us consider visuals and sounds (Figure \ref{natural}). To each visual object (a simple shape as shown in Figure \ref{basis}), we associate a coefficient. 
The superposition of simple visual shapes with their coefficients constitute the filtered visual images. Visual coefficients tell us {\em where} an image is within a space and which size it has.\footnote{We could also talk of `dimension' here, however specifying that this is about the size of each ket in the space, not about the number of elements in a basis.} Thus, the objects represent {\em which element} (e.g., a circle, a point, a straight line) do we need, and the arrows tell us {\em where, how big} and possibly also {\em how many} simple shapes of a specific kind do we need (e.g., where, how big, and how many circles do we need).  Let us now consider sounds, where the objects are the simple musical sequences (obtained as sonified filtered images). To each simple musical sequence we can associate a coefficient. Sound coefficients tell us {\em when} a musical sequence is played, how loud/long it is. Coefficients give a structure to the set of images or sounds, giving them an organization. Each (visual) coefficient represents the action of taking an image from the set of all kets and putting it into a point of the space with a certain size. Thus, visual coefficients, whose action can be indicated via an arrow, act as positioning and resizing. Similarly, auditive arrows act as positioning in time and resizing in loudness.
We have a correspondence between coefficients about position in tridimensional space and size, and coefficients about position in temporal space and loudness. The transformation (indicated here as the transform $T$) from the domain of visuals to the domain of sounds associates images with musical sequences, and visual coefficients with sound coefficients. If we consider two different sonification strategies, let they be $T$ and $T'$, for example based on slightly different choices of vectors, we can define  transformations from musical pieces obtained via $T$ and musical pieces obtained via $T'$. 
This is equivalent to compare two musical pieces, one obtained via a given technique, and the other obtained via another one; see Figure \ref{natural}.

While decomposing the original image, we can distinguish between a pattern (e.g., if we consider a picture of a fish, the fish's ``skeleton'') and an envelope (e.g., the overall shape of the considered fish). The envelope is represented by the $0,...,j$-th coefficients, while the pattern by the other
 $j+1, ..., N$-th coefficients. If we start from $0$, formally the last index should be $N-1$ for the $N$-th coefficient.
We have:
$|\phi_i^v\rangle=\left\{ ( |\phi_1^v\rangle, |\phi_2^v\rangle, ..., |\phi_j^v\rangle), (|\phi_{(j+1)}^v\rangle, |\phi_{(j+2)}^v\rangle, ..., |\phi_N^v\rangle) \right\}$.
The same applies to sound:
$|\phi_i^s\rangle=\left\{ ( |\phi_1^v\rangle, |\phi_2^s\rangle, ..., |\phi_j^s\rangle), (|\phi_{(j+1)}^s\rangle, |\phi_{(j+2)}^s\rangle, ..., |\phi_N^s\rangle) \right\}.$
If we consider together the two notations, we can run a double index:
$|\phi_i^{\mu}\rangle$, where $\mu = 1$ stands for $v$ (visuals), and $\mu = 2$ stands for $s$ (sounds).

After having defined our bases, we can talk about the measurement idea. In quantum mechanics, a measurement of a state implies a destruction of part of the initial information (destructive measure). The reason is that a quantum state, described by a wave function, is usually a superposition of eigenstates, that are reduced to a single state after the measurement. The information that is preserved depends on the measurement apparatus. Intuitively, let us consider two states, with spin-up and spin-down respectively, each one with a probability density of $\frac{1}{\sqrt{2}}$, and thus, probability $\frac{1}{2}$. If the experimental apparatus is prepared to detect a specific spin, say, spin-down, after the measurement the spin-up component will be deleted, and the wave function will {\em collapse} to spin-down. Thus, after the measurement the state is only constituted by a spin-down with probability $1$, and all future measurements with the same measurement apparatus will give spin-down with probability $1$. A practical application of this idea is shown by the Stern-Gerlach experiment \citep{cohen}.

Here, we are not dealing with spin but we can borrow the essential idea about quantum measurement.
In the framework we proposed, our ideal ``measurement instrument'' is constituted of a selection of visual (and sound) kets. For example, if we have the image of a fish and a few kets' superpositions giving an ellipse, a triangle, and a circle, we will have lost all the information about the original fish image except the contour that is schematized via an ellipse, the tail with a triangle, and the shape of the eye via a circle, and maybe the fins with some more triangles. In this sense, the image has been filtered via the an ellipse, some triangles, and a circle. If we make other measurements of the same filtered fish, we will get the same image, because it has already been filtered via that filter. The circle, the ellipse, and the triangle(s) will have a specific coefficient that describe in which spatial point the simple shape is put with which size and so on.
For this reason, measuring in this framework is like a filtering process; more precisely, expressing something in terms of (a finite number of) basis coefficients. If we want to sonify our filtered fish, we have first to use a ``translation dictionary'' (that verifies gestural similarity condition) to see how to sonify each of the simple shapes used, and then we have to superimpose these musical sequences. Gestural similarity condition also applies to the position and size of each one of the selected visual kets. Summarizing, we measure an image, and then we take its transform into sound. 


The first rough approximations of the original shape are just in terms of essential lines.
This is the envelope part:
\begin{equation}
T\left( \sum_{i=0}^j|\phi_i^v\rangle \right) = \sum_{i=0}^j|\phi_i^s\rangle
\end{equation}
However, if we map only the patterns, we'll just have the terms:
\begin{equation}
T\left( \sum_{i=j+1}^N|\phi_i^v\rangle \right) = \sum_{i=j+1}^N|\phi_i^s\rangle
\end{equation}
This requires that matrix $T$ has some block form, so we exclude index permutations.
We can get two partial sonifications: a sonification of the envelope and a sonification of the pattern.
The complete musical rendition would require summing up (superposing) these terms:
\begin{equation}
T\left( \sum_{i=j+1}^N|\phi_i^v\rangle \right) + T\left( \sum_{i=0}^j|\phi_i^v\rangle \right)
\end{equation}
 
 

 
 \begin{equation}
 \begin{split}
 |\psi^v\rangle \sim \sum_{i=0}^Nb_i|\phi_i^v\rangle = & \sum_{i=0}^j b_i|\phi_i^v\rangle + \sum_{i=j+1}^Nb_i|\phi_i^v\rangle
  \Rightarrow  T\left(\sum_{i=0}^j b_i|\phi_i^v\rangle\right) + T\left(\sum_{i=j+1}^Nb_i|\phi_i^v\rangle \right)
 \\ = & T\left(\sum_{i=0}^j b_i|\phi_i^v\rangle + \sum_{i=j+1}^Nb_i|\phi_i^v\rangle \right)
 = T\left( \sum_{i=0}^Nb_i|\phi_i^v\rangle \right)\sim T(|\psi^v\rangle)\sim |\psi^s\rangle
 \end{split}
 \end{equation}
 Thus: $T|\psi^v\rangle \sim |\psi^s\rangle$, and $T:[\textbf{v}]\rightarrow{[\textbf{s}]}$, where the new notation with square brackets is meant to highlight the basis change. While the original image can be seen as the product of envelope and pattern, the final result after filtering and base decomposition (that is, our selection of visual kets to decompose a given image) may be seen as the {\em sum} of coefficients, as terms within a matrix. The sum is direct because there is no overlap between the meaning expressed by coefficients.
Thus, the external product between overall shape and pattern is transformed into a direct sum:
 $(\mbox{overall shape}\otimes\mbox{patterns})\Rightarrow ((0\rightarrow j\mbox{ coefficients }) \oplus ((j+1)\rightarrow N \mbox{ coefficients }))$.
The external product becomes a direct sum of visual kets combined via suitable coefficients. In fact, we can represent the basis as follows: $|\phi_0^v\rangle,...,|\phi_j^v\rangle\bigoplus|\phi_{j+1}^v\rangle ,...|\phi_{N}^v\rangle$
This basis can be represented as a matrix. Sonification operator $T$ acts on coefficients from $0$ to $j$ and from $j+1$ to $N$, getting  $|\phi_0^s\rangle,...,|\phi_j^s\rangle\bigoplus|\phi_{j+1}^s\rangle ,..., |\phi_{N}^s\rangle$.

Summarizing, we can decompose the initial image via a basis of visual vectors: $|\psi^v\rangle\sim\sum_i^Nb_i|\phi_i^v\rangle,$
and then, we can map this decomposition into sounds: $ T(F|\psi\rangle)=T\left( \sum_i^Nb_i|\phi_i^v \right)=\sum_i^Na_i|\phi_i^s\rangle\sim|\psi^s\rangle,$
provided that each element of the visual basis corresponds to one element of the sound basis.


If the final music is refined (operator $R$; e.g., well-orchestrated, with motives and chords completed, and so on), we have: $R(T(F|\psi^v\rangle)) = |\psi^s\rangle, $ with $|\psi^s\rangle$ a complete, ready-to-perform piece of music.
Thus, we have Equation \ref{RTF}; see on the right a diagrammatic rendition of the same equation. Summarizing, a possible protocol to make a complete musical piece out of an image is is filter--transform--refine, or F-T-R. Eq. \ref{RTF} contains the operators in inverse order because they are applied from right to left, according to the (mathematical) composition convention also used in category theory while composing arrows. See also Figure \ref{RTF_image}.
\begin{equation}\label{RTF}
\begin{diagram}
|\psi^s\rangle = RTF|\psi^v\rangle,\,\,\,\,\,|\psi^v\rangle & \rTo^{RTF} & |\psi^s\rangle
\end{diagram}
\end{equation}
The refinement action is not unique; it is subject to a choice --- also the basis is a choice. 
$R$ maps the superposition of sound vector basis within a complete piece. If we consider all possible (not-too-free) choices for $R$, we grasp an intuitive idea of the equivalence class of musical pieces
that can be composed via a refinement of the same starting material. Differences between these pieces are less evident if we neglect $R$ and we just consider several, possible bases for the first step, visual decomposition. These pieces are the sonification(s) {\em up to an isomorphism} of the initial shape. In fact, if we evaluate the distance between different sonification/musical renditions of the same simplified shape, we expect to find a value $\eta$ (eq. \ref{distance}) that is different from the distance between the sonification of a shape and the sonification a completely different shape. This leaves room for experimental investigations of shape differences and values of musical distances. These last ones can be evaluated through matrix distance evaluations, as done in \citealt{mannone_compagno}. Equation \ref{distance} shows the information we would like to obtain from a suitable distance criterion, where $\alpha$ and $\beta$ are two visual decompositions. These decompositions are {\em minimal} for Gestalt: in both cases there is only the minimum number of simple shapes that, combined together, make the image recognizable. The quantity $\epsilon$ is a positive, real number. If we apply $T$, we get the criterion for musical distances, see the second part of eq. \ref{distance}, with $\eta$ a real and positive number (we do not use $\delta$ here to avoid confusion with the Kronecker delta). 
\begin{equation}\label{distance}
\mathcal{D}\left( \sum_ib_i^{\alpha}|\phi_i^{\alpha,v}\rangle,  \sum_ib_i^{\beta}|\phi_i^{\beta,v}\rangle\right)<\epsilon,\,\,\,\,
\mathcal{D}\left( T\left(\sum_ib_i^{\alpha}|\phi_i^{\alpha,v}\rangle\right), T\left( \sum_ib_i^{\beta}|\phi_i^{\beta,v}\rangle\right)\right)<\eta
\end{equation}
The distance $\mathcal{D}$ defined in eq. \ref{distance} can be formulated in terms of cosine similarity, drawing a parallel with studies on (musical) tonal attraction via quantum formalism \citep{graben}. With cosine similarity, the distance $\mathcal{D}$ is $0$ if the two basis decompositions (and their sonifications) are completely unrelated, and $1$ when they are identical. In this framework, $\epsilon$ and $\eta$ should be positive real numbers much smaller than $1$.

Eq. \ref{distance} may suggest that $T$ is a continuous function. However, the general idea is that different, but equivalent, base decompositions of an image should have a small distance between them, and their transform $T$ should keep such a distance small. It means that, within the variety, the listener should be able to connect equivalent filtered images with equivalent sonifications, where equivalence is evaluated within the limits of $\epsilon$ and $\eta$.
The value of $\epsilon$ is characteristic of the given original image, and it can be evaluated experimentally via the cognitive experiment described in the following section. Ambiguity leads to higher values of $\epsilon$: that is, if two initial images are quite close or they show some visual ambiguities, also their decompositions will do.

Coefficients $a_i$, $b_i$ indicate where the corresponding sound or visual ket should be put. In fact, the same kets can create different images if recombined in different positions or with different relative size. To simplify our analysis here, however, we will just consider $a_i$ and $b_i$ as scalars, and thus, as real numbers.
\begin{equation}
\langle \phi_i^v|\phi_j^v\rangle = \delta_{i,j},
\end{equation}
where $\delta_{i,j}$ is a Kronecker delta, with values $0$ or $1$.
With vectors of different basis:
\begin{equation}
\langle \phi_i^v|\phi_{j}^{'v}\rangle = \epsilon_{i,j},
\end{equation} 
where $\epsilon_{i,j}$ is a positive number comprised between $0$ and $1$ that indicates the distance between corresponding elements of the two bases.
Thus, $\delta_{i,j}\in\left\{ 0,1 \right\}$ and $\epsilon_{i,j}\in[0,1]$.
For the same initial image, two different but {\em equivalent} decompositions will either have the same visual kets with slightly different coefficients, or also different visual kets but with a negligible difference. In fact, this $\epsilon_{i,j}$ is the piecewise corresponding to the overall $\epsilon$ as described above, that is supposed to be much smaller than the distance between the decompositions of different images. The same applies to sonification, and thus to the selected kets with index $s$.
If we compare a selection of (sonified) visual kets, that is, sound kets of a sonified visual basis and sound kets of another sound basis:
\begin{equation}
(T^{\ast}\langle\phi_i^{v}|)|\phi_{j}^{'s}\rangle=\langle \phi_i^s|\phi_{j}^{s}\rangle = \eta_{i,j},
\end{equation} 
where $T^{\ast}$ is the adjoint of $T$ acting on the bra vector.
A comparison between the original image and the filtered image will give an idea of the accuracy of filtering, and of the amount of information we can ignore to still leave the original shape easy to recognize. In the same way, once we have a refined musical piece, we can compare it with the superposition of sound kets, and we can evaluate the amount of personalization the composer put, and what are the notes and musical sequences he or she probably started from. 
  

Filtering operation, represented by operator $F$, can be compared with quantum measure, that is a destructive measure. This would be one more reason to justify the use of a quantum-like formalism adopted here.

Gestural similarity conjecture (see Section \ref{introduction}) applies to the mapping of each element of the visual basis to each element of the sound basis (see diagram \eqref{gest_diagram}), where $\mathcal{G}$ is the common gestural generator that produces an element of the visual basis and its corresponding element of the sound basis, and $\forall i=0,...,N$.
\begin{equation}\label{gest_diagram}
\begin{diagram}
  & & \mathcal{G}  & &
\\ & \ldTo^{generates} & & \rdTo^{generates}
\\ |\phi_i^v\rangle & & \rTo^T & & |\phi_i^s\rangle
\end{diagram}
\end{equation}


This approach can include temporal evolution. In music, $|\phi_i^s\rangle$ already contains time: in fact, each simple shape is seen as a collection of sound events in time. Thus, there is a time where each sound sequence/melody/superposition of short melody lives. If the initial shape is changing over time, we will have a collection of sound sequences, each one for each time instant. To make sonification practical and possible, we can consider snapshots of the changing image at given time intervals. Thus, we can obtain a finite sequence of sound sequences per shot. This would correspond to the sonification of a chronophotography of the changing original image.
Thus, $|\phi_i^s(t)\rangle$ is the sound sequence at each step $t$, where $t$ indicates when the shot has been taken. This means that the transformation of $|\phi_i^v(0)\rangle\rightarrow|\phi_i^v(T)\rangle$ of the image, from time $t=0$ to time $t=T$, correspond to a sound transformation $|\phi_i^s(0)\rangle\rightarrow|\phi_i^s(T)\rangle$.
Temporal operator $\mathcal{T}$ acts as: $$ \mathcal{T} |\phi_i^s(0)\rangle=|\phi_i^s(t)\rangle. $$ More precisely, we can split the action of temporal operator between sounds and visuals: $\mathcal{T}_s\otimes\mathcal{T}_v$, where $\mathcal{T}$ can be freely compared with a tensor (or monoidal) functor. Temporarily neglecting coefficients, we have:
$$ T((\mathcal{T}_v\otimes\mathcal{T}_s)|\phi_i^v\rangle) = T(\mathcal{T}_v|\phi_i^v\rangle\otimes\mathcal{T}_s(1)) =T(\mathcal{T}_v|\phi_i^v\rangle)=T(\mathcal{T}_v\otimes\mathcal{T}_s|\phi_i^v\rangle)=\mathcal{T}_sT(\mathcal{T}_v|\phi_i^v\rangle)=\mathcal{T}_v|\phi_i^s\rangle, $$
and
$$ T(\mathcal{T}|\phi_i^s\rangle)=T((\mathcal{T}_v\otimes\mathcal{T}_s)|\psi_i^s\rangle)=T(\mathcal{T}_v(1)\otimes\mathcal{T}_s|\phi_i^s\rangle)=T(\mathcal{T}_s|\phi_i^s\rangle),$$
where $I$ indicates a unit vector. 
But $T$ applied to sounds acts as the identity, thus $T(\mathcal{T}_s|\phi_i^s\rangle)=\mathcal{T}_s|\phi_i^s\rangle$.

A possible, computational goal of this approach could be the development of an app that takes as input an image, and sonifies it in real time. If the app allows a user manipulation of the image, with, for example, stretching and scaling (and any other homotopic transformation), or cutting/gluing/repeating (and any other non-homotopic transformation), the music should change in real time. This conjecture may lead to a conservation theorem.
Because conservation theorems in physics are deeply linked to the concept of symmetry as a ``change without change'' \citep{wilczek}, and symmetry is traditionally associated with beauty in nature and in the arts \citep{weyl}, we can hope that seeking conservation and symmetry principles in music and visual arts could allow us to develop stronger tools to analyzing existing artworks and creating new ones. 

\begin{figure}
\centering
\includegraphics[width=6cm]{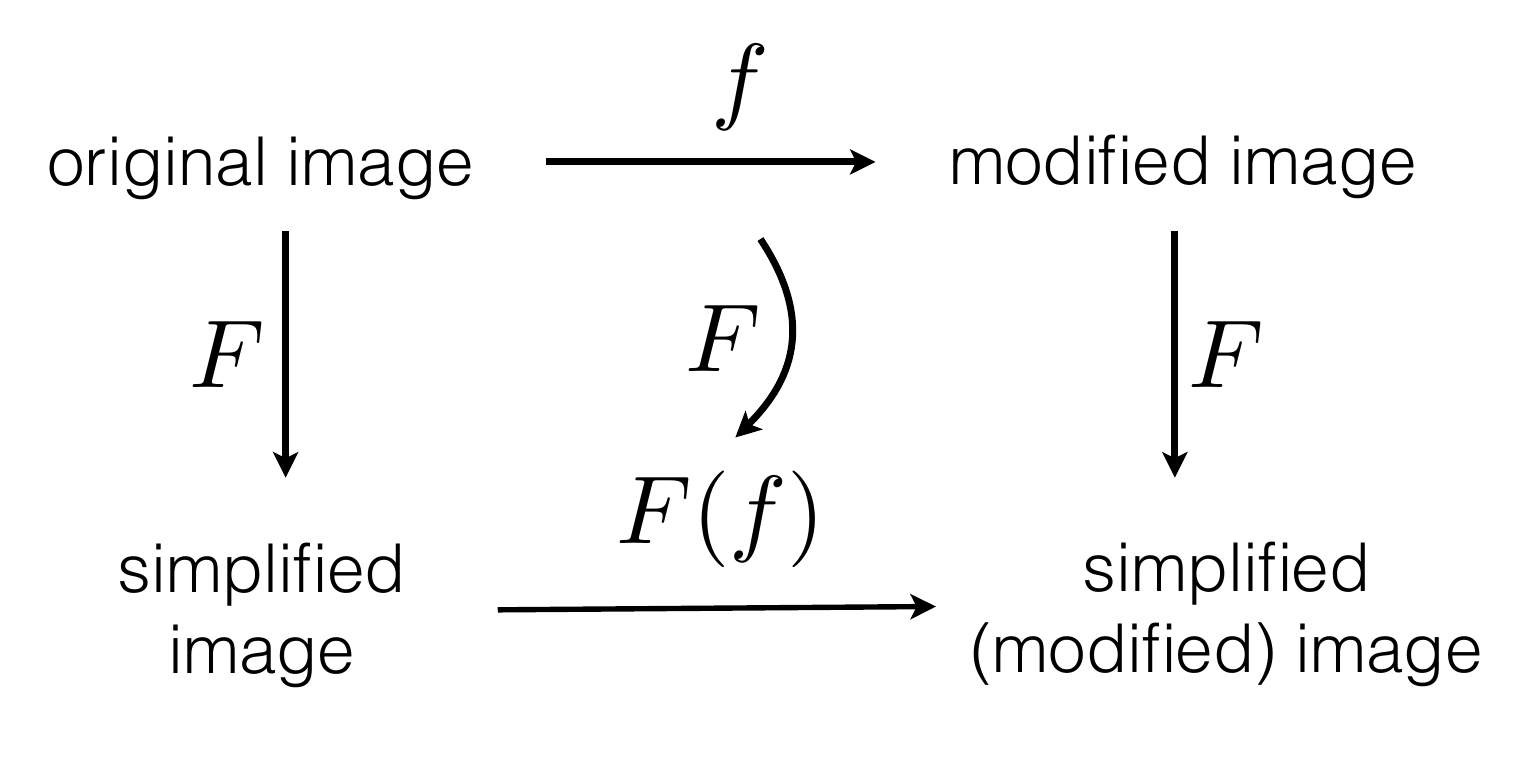}
\caption{An image is simplified via $F$. The image can also be slightly modified via a sort of multiple transformations indicated via $f$.
We can imagine to simplify the original image via filter $F$, and to simplify the modified image via $F$. The transformation between the two filtered images can be seen as the composition of $f$ with $F$. We can also think of a class of transformations $f$ mapped onto a class of transformations $F(f)$, as well as a class of original images and a class of modified images mapped onto a class of simplified images and a class of simplified modified images, respectively.
}
\label{method1}
\end{figure}

\begin{figure}
\centering
\includegraphics[width=8cm]{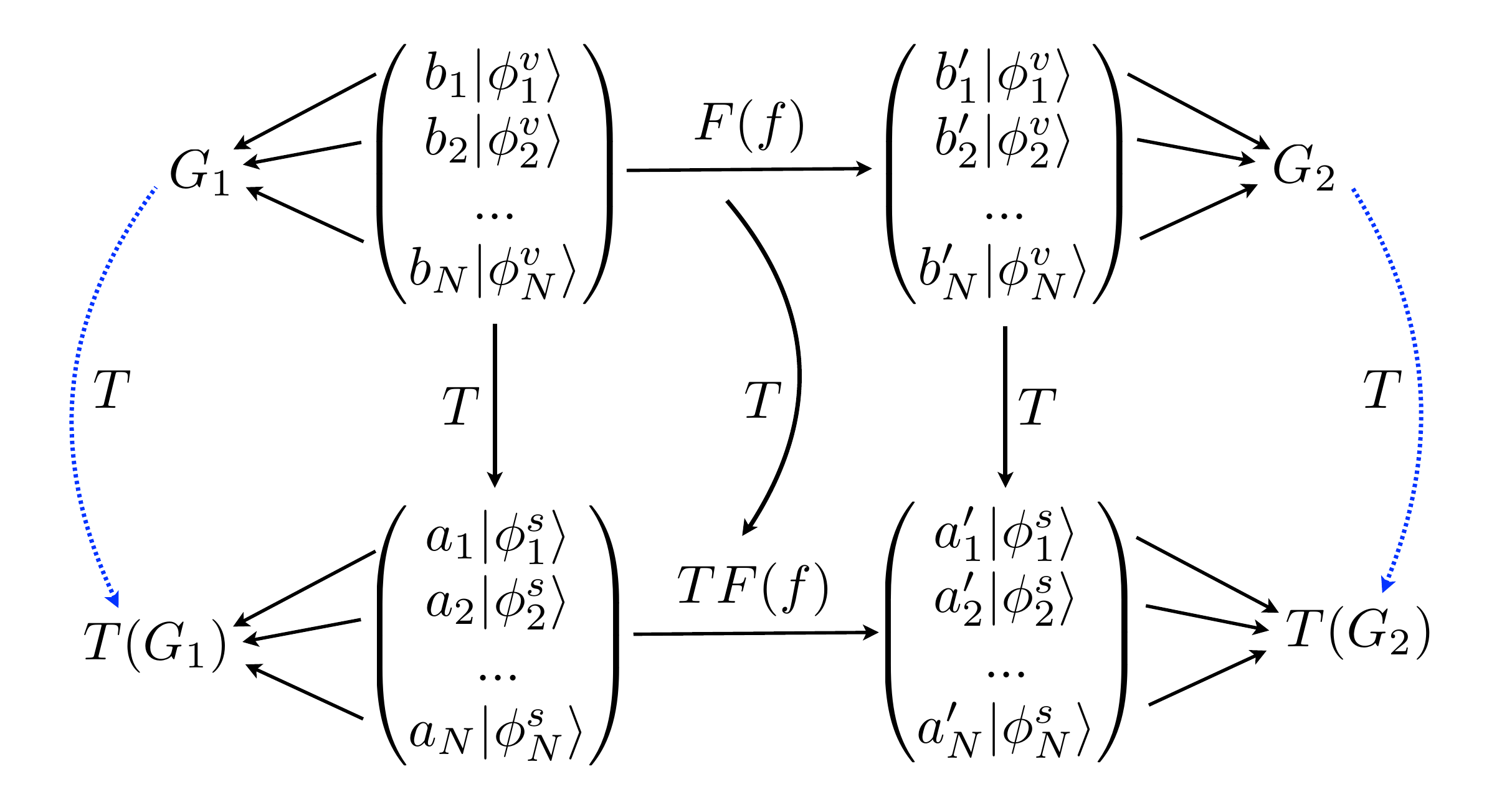} 
\caption{Each simplified image can be formally described via a superposition of visual kets. In the figure, we schematically represent the selected kets and the coefficients via an overall column vector. Transform $T$ also maps basis decompositions (simple shapes with visual coefficients giving them a structure in visual space) in the domain of visuals into basis decompositions in the domain of sounds (musical sequences with coefficients giving them a structure in musical space). If we consider $f$ as a class of transformations $f$ as described in Figure \ref{method1}, thus we can include in our description also the (class of) deformations from a simplified image into another simplified image as $F(f)$, and the comparison between a sonified simplified image with the sonified deformed shape via $TF(f)$. Finally, the symbols $G_1$ and $G_2$ are introduced here to represent the overall Gestalt of the first image and the second image, respectively.}
\label{method2}
\end{figure}



\begin{figure}
\centering
\includegraphics[width=8cm]{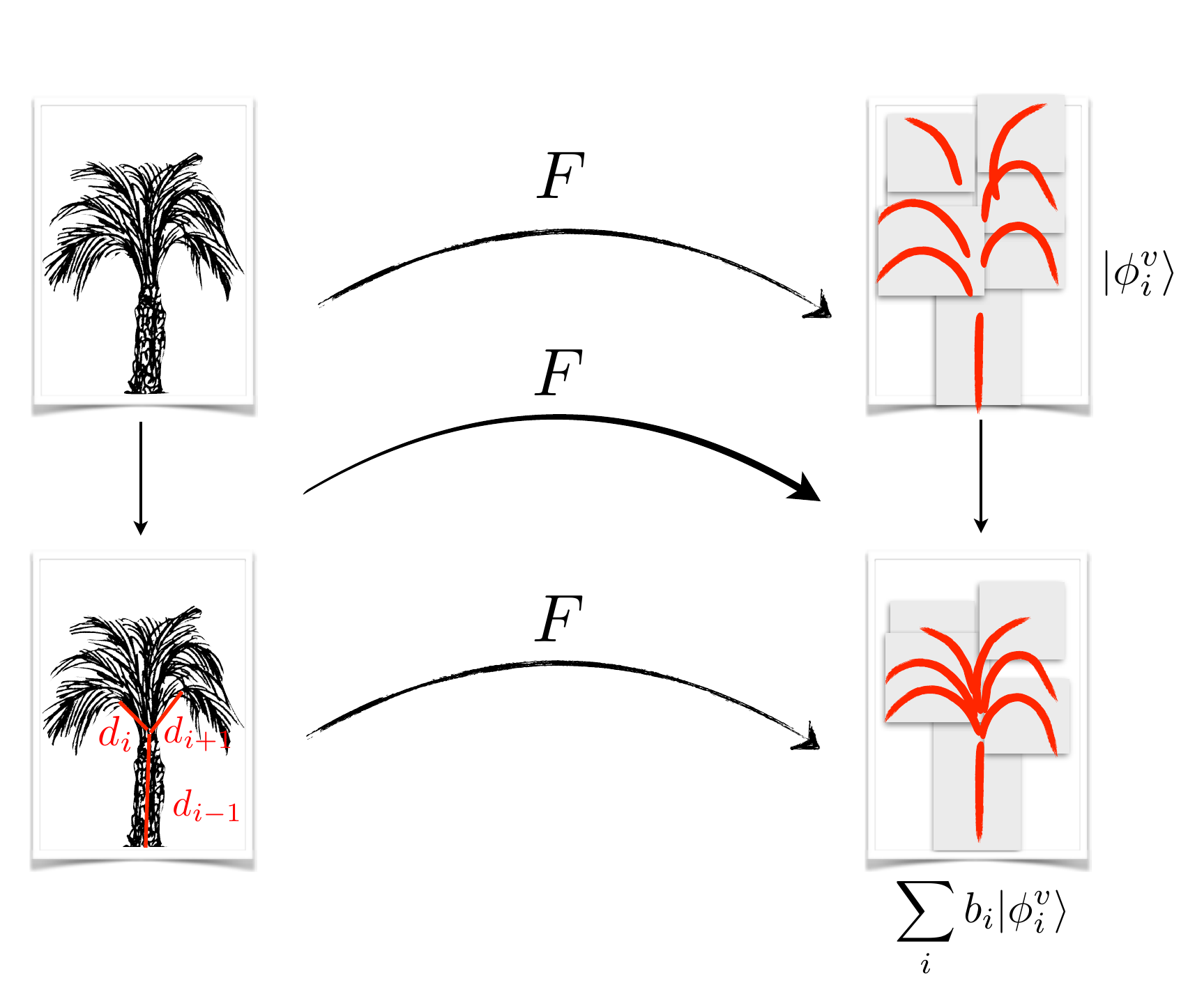} 
\caption{
$F$ maps original images onto their decomposition in terms of `covering' charts made of simple shapes. If we define distances $d_i$ between points of the original image as connected paths, we can map these distances via $F$ to give a structure to the charts.
}
\label{charts}
\end{figure}

Figures \ref{method1} and \ref{method2} illustrate the described procedure.
We can wonder if a visual Gestalt can be translated into a sound Gestalt: formally, this means asking if the blue, dotted arrows in Figure \ref{method2} are preserving the structure of Gestalt from one domain to the other.  Is it possible/easier to translate a visual Gestalt into a musical Gestalt? The correspondence between visual kets and sound kets that were supposed to be joined by gestural similarity has been recently verified \citep{collins, mannone_collins}. However, the correspondence of the initial Gestalt and of the final, overall Gestalt, that involve superpositions of kets with appropriate coefficients, has yet to be proven. Section \ref{experiment} proposes an idea for such an experiment. If successful, this would extend the action of operator $T$ to Gestalt(s), making the overall diagram commutative.
Finally, let us focus on filtering $F$.
We can associate to each pair of points of the original visual image their distances $d_i$ between them in terms of connected path and composition of connected paths.
Thus, each visual image can be enriched with all distances, see Figure \ref{charts}. The image obtained via $F$ can be thought of as a covering, ``like a triangulation'' with charts (made by simple shapes) on the original image. Thus, $F$ maps distances $d_i$ between points of the original image onto distances between charts; this helps build a structure out of the charts, constituting a sort of atlas. In this way, we could in principle connect our formalism with hints from differential geometry.    

\begin{figure}
\centering
\includegraphics[width=10cm]{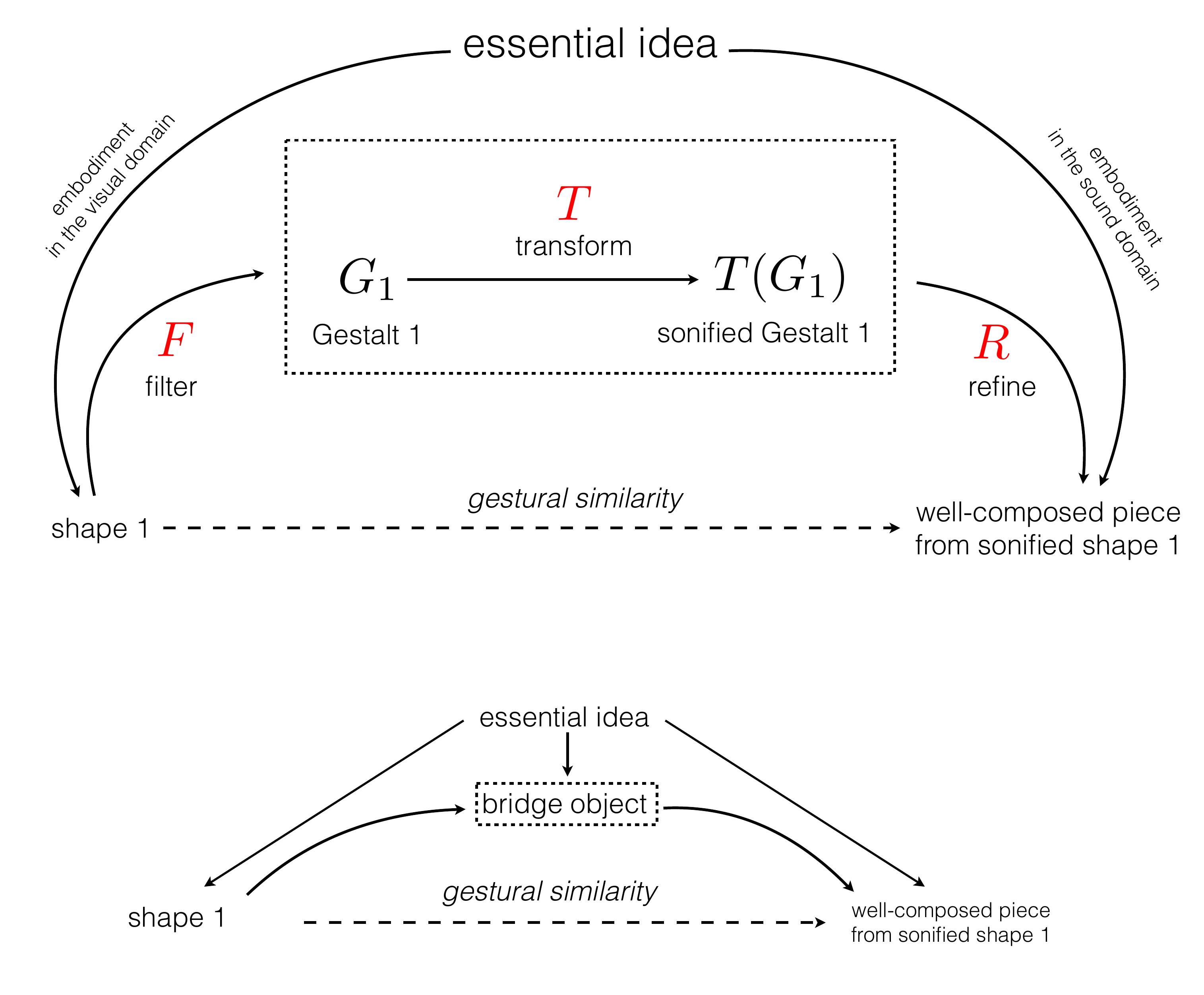}
\caption{The text in the dotted box is freely inspired of Caramello's bridge \citep{caramello2}, that has been recently joined with gestural similarity/essential idea \citep{mannone_favali}. The red letters highlight the discussed $RTF$ procedure, that constitutes the overall bridge.}
\label{bridge}
\end{figure}

Overall, our approach can be compared with the definition of ``connecting bridges''\footnote{According to {caramello2}, ``the importance of `bridges' between different areas lies in the fact that they make it possible to transfer knowledge and methods between the areas, so that problems formulated in the language of one field can be tackled (and possibly solved) using techniques from a different field, and results in one area can be appropriately transferred to results in another.'' More specifically, Caramello refers to topoi as bridges to unify different areas of mathematics; however, the bridge metaphor can be generalized.} by Olivia \cite{caramello2}. In fact, a possible Gestalt translation reminds us of the ``bridge object'' discussed in Caramello's theory and recently applied within a gestural similarity framework \citep{mannone_favali}.
We can start from visuals and analyze the Gestalt extracting a core of Gestalt information as an independent object --- that, in Caramello's framework, would be metaphorically called a bridge object --- and this core could be applied within the musical domain, constituting a basic composition whose Gestalt is, hopefully, the musical rendition of the initial one. The lower part of Figure \ref{bridge} represents the relationship between the approach with an essential idea that can be embodied within the visual or the sonic domain, and a bridge object, for example a mathematical formalization of some shape properties, that joins visual and sounds. The upper part of Figure \ref{bridge} shows that our approach, with a possible Gestalt translation, can be contextualized within the bridge object, where the entire bridge is constituted by the composition RTF, with filter $F$, transform $T$, and refinement $R$. We start from the domain of the continuous, with the original image; then, we take a discrete version of it with the filter, and, finally, with the refinement step, we get back to the continuous. A sonified filtered image (the result of $TF$) can be thought of as a discrete superposition of sonified visual kets, that is, specific sound kets. If a composer starts from it and he or she adds more notes and refines passages, the result of the $TF$ sequence is turned into a complete and aesthetically satisfactory musical piece. This is the refinement step ($R$), and it can be thought of as an interpolation.

If the work is done well, we might expect that the final result will remind the listener of some of the essential characteristics of the initial shape. Thus, the gestural similarity between the very initial shape and the very final shape will be verified, and the overall diagram will be commutative. 
Of course, if we reverse all the arrows, we can get the converse, that is, visual shapes from music \citep{mannone_libro}. We can wonder if we can go from image to sound and back, and what is the entity of the information loss. This also depends on the amount of information encoded in the coefficients during the filtering step.
If we sonify and image and then we want to reconstruct the original one, we do expect a partial and approximate reconstruction, because of the (two) filtering passages. We expect to keep -- and we can appropriately select parameters to keep it -- the general recognizability of the original image, or, at least, of the equivalence class to which the image belongs. Time/loudness/pitch domain for music and (visual) space domain for images are independent, but they can be related via the reading-through-time of the (visual) space elements. We can consider a global time, corresponding to the overall (visual) space, with some local times in it, flowing within each small figure or simple visual element. This can be compared with musical fragments performed by each musician during a whole symphony. Also, the visual-sound correspondence is more robust if we limit ourselves to drawing, or if we simplify complex visuals with collections of lines. This is consistent with some recent experimental approaches \citep{collins, mannone_collins,thoret}.

\subsection{Gestalt and topology}

\begin{figure}
\centering
\includegraphics[width=6cm]{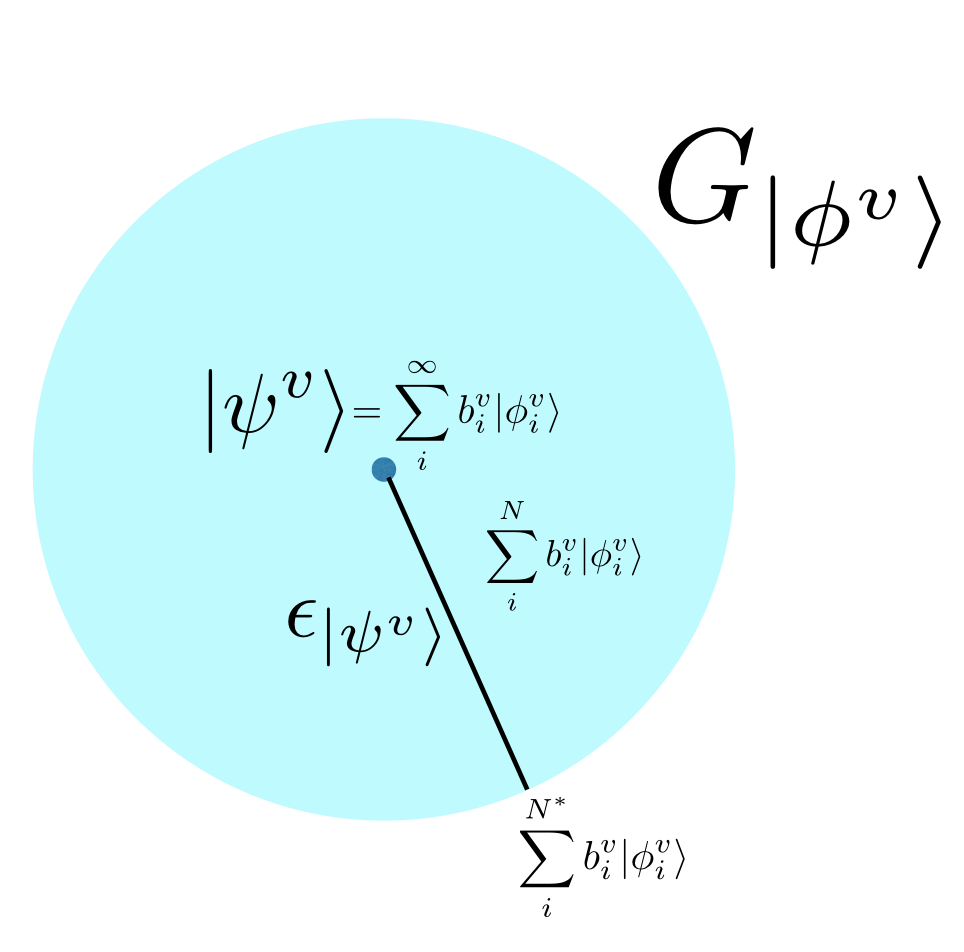}
\caption{A topological Gestalt neighborhood $G$ of a visual state $|\psi^v\rangle$. The closer we are to the center, the higher the number of base elements is, and the closer the decomposition is to the original state. Conversely, the farther we are from the center, and the more coarse-grained the approximation is. With a number of elements less than $N^{\ast}$, the original image is no longer recognizable. The same idea applies to sound states as well. Bistable states can be thought of as two apparently superposed circles, as the two faces of a cylinder with a small height. We cannot look at both faces at the same time. In the Figure, the choice of color is arbitrary.}
\label{ball}
\end{figure}

Through experiments, we can find out which, and  how many visual kets are needed to recognize a given image; that is, what is the smaller index $N^{\ast}$ necessary to recognize an image or part of it. A similar criterion can be adopted to find the minimum amount of base elements to recognize a melody or a generic musical piece. If we reduce musical pieces to their skeletal structure, a number of same-genre compositions are expected to be indistinguishable, as happens with Schenkerian analysis. An extension of our quantum-like approach can involve the topological definition of neighborhood. We can define a Gestalt-neighborhood of a state $|\psi^v\rangle$ for visual (the same idea applies to sound states), where the center is the original image (that can be perfectly approximated with infinite coefficients), and where the most external layer is represented by the decomposition with the least amount of elements, $N^{\ast}$ elements. A decomposition with a number $n\in\mathbb{N}$, with $n<N^{\ast}$ of elements, will produce ambiguities with other images; see Figure \ref{ball}.

\subsection{Examples of visual and sound kets}\label{examples}

\begin{figure}
\centering
\includegraphics[width=13cm]{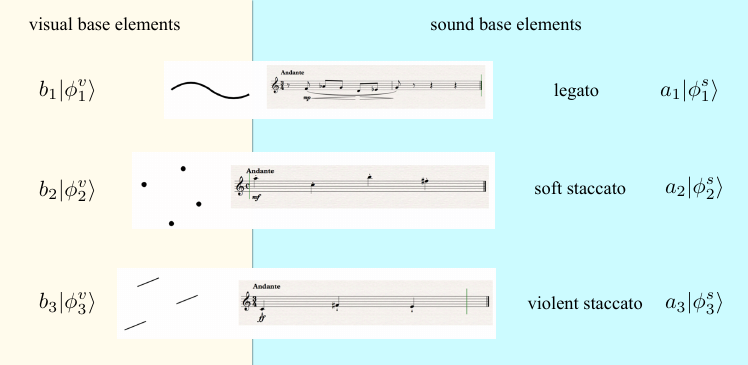}
\caption{Let us consider three visual kets: let $|\phi_1^v\rangle$ be a dot, $|\phi_2^v\rangle$ a short segment, and $|\phi_3^v\rangle$ an arch. A very simple superposition of visual kets of the same kind gives, in this example, a curved line, a sequence of dots, and a sequence of angled segments, respectively. Because, in this case, the superpositions concern elements of the same kind, we can see them not as sums but just as the results of the action of coefficient $b_1$ for the first visual ket (i.e., ``connect two arches to get a curved line,'' that is, ``repeat the gesture of creating short curved lines in the described way''), coefficient $b_2$ for the second visual ket (i.e., ``repeat the gesture of creating dots with this distance''), and coefficient $b_3$ for the third visual ket (i.e., ``create the segment with this inclination as the result of a nervous gesture and repeat it with this distance''). These images (and many others), and the corresponding sounds (and many others), have been used in an experiment \citep{collins, mannone_collins} to measure the degree of the assessed gestural similarity. The results confirmed the theoretical expectations: in this example, a curved line is associated with a {\em legato} melody independently by the speed or the pitch range; a sequence of dots is associated with {\em staccato} notes independently by their pitch range or their speed or loudness, and, finally, a sequence of angled straight segments is associated with a {\em staccatissimo}, a ``violent'' {\em staccato} sequence of notes, independently by their pitch range \citep{collins, mannone_collins}. The main parameters involved in gestural similarity in the analysis run so far are pitch contours -- not within a specific pitch range -- and articulations.   
There is a correspondence between a few elements of the visual base, and some elements of the sound base.
In the pseudocode of Section \ref{pseudocode}, pitches of the sound sequence can be customized. }
\label{basis}
\end{figure}

The proposed method can be used to obtain concrete musical sequences from a given image. The preliminary condition is the definition of a function from image to music, and the use of a dictionary of simple shapes and simple sound sequences connected via gestural similarity.
An appropriate selection of suitable elementary visual kets and of corresponding sound kets should involve some experimental validation. For example, we can start from some hypothesis, based on gestural similarity, and we can check it via an experiment. We can give people a picture of a clear form, for instance a palm, and a set of possible simple shapes, with the possibility to repeat a shape several times. Which shapes will be more often used for approximately sketching or tiling the given form? A similar experiment could involve an excerpt of a musical piece and a collection of short sound sequences. Which short sound sequences will be more often used to approximately sketch and represent the given musical piece? Another experiment would, of course, involve the correspondence between visual kets and sound kets.
Initial validations with simple visuals and simple musical sequences has been done \citep{collins, mannone_collins}; see Figure \ref{basis} for some of the visual and sound examples used in that cognition experiment. In Figure \ref{basis}, the visual kets could be constituted by a dot, an arch, and a segment. Their combination (with repetition, angle selection, position in the space, formally expressed by coefficients $b_i$) generates a curved line, a collection of dots, and a sequence of angled segments, respectively. Each of these combinations appears as being generated by a specific gesture, and the same gesture can be thought of as being at the origin of the sound sequences on the right of the image. A good correspondence between these simple visual sequences and musical sequences has been proved \citep{collins, mannone_collins}.

Summarizing, we need both visual and sound kets with their coefficients. Coefficients for the visual kets elements characterize position and size, and their information will be translated into onset, duration, pitch range  -- and eventually also loudness, for tridimensional images \citep{mannone_libro} -- values for the coefficients of sound kets. 
Also, we can use a (mathematical) composition of visual kets: e.g., a dotted spiral can be considered as the composition of a spiral form (this is an example of what we call here ``envelope'') and of dots (``pattern''). Musically, we will have an ascending-descending (along the spiral) sequence of notes, with {\em staccato} articulation. The staccato articulation corresponds to the isolated dots, while the spiral form can be translated into an equivalent pitch movement in the space of pitch. In fact, aside from time which could be considered as overall flowing from the left to the right for an image, the small images a complex image is made of (the small images a complex image can be decomposed into) present a local time that flows on them. For example, if we see a straight line and a spiral, we can depict a melodic line following the straight line, and another melodic line following the profile of the spiral. This is intrinsically different from a simple scanning of an image from left to the right, as it can be the sound rendition of the spiral within the spectrogram of Figure \ref{Windowlicker} --- where, as an additional difference, time-frequency space is used rather than time-pitch space. We can also think of representations including a spectral dimension.

\section{A pseudocode and computational developments}\label{pseudocode}

In this section, we propose a possible computational application of the theoretical method described in Section \ref{mathematical details}. This procedure can be contextualized within the framework of algorithmic composition. The use of a set of rules to compose music is not new,
and the use of computer has quite a long tradition as well \citep{xenakis}.
In algorithmic music, part of the compositional process is delegated to the machine, but the human creativity relies on the choice itself of the process and of the specific parameters that can dramatically change the final outcome. Also, in our method, images and their parameter choices, as well as the number of base elements to schematize them, play a relevant role.

\begin{figure}[h]
	\centering
		\includegraphics[width=0.3\linewidth]{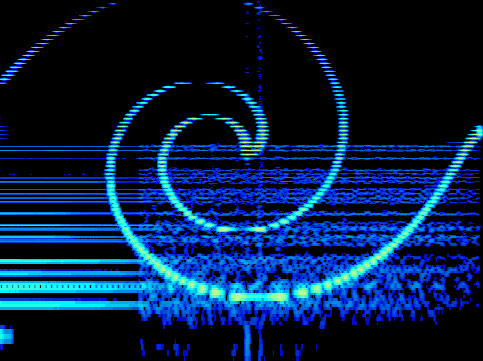}
		\caption{A drawing within a spectrogram. Image from \cite{niinisalo}.}
		\label{Windowlicker}
		
	\includegraphics[width=0.6\linewidth]{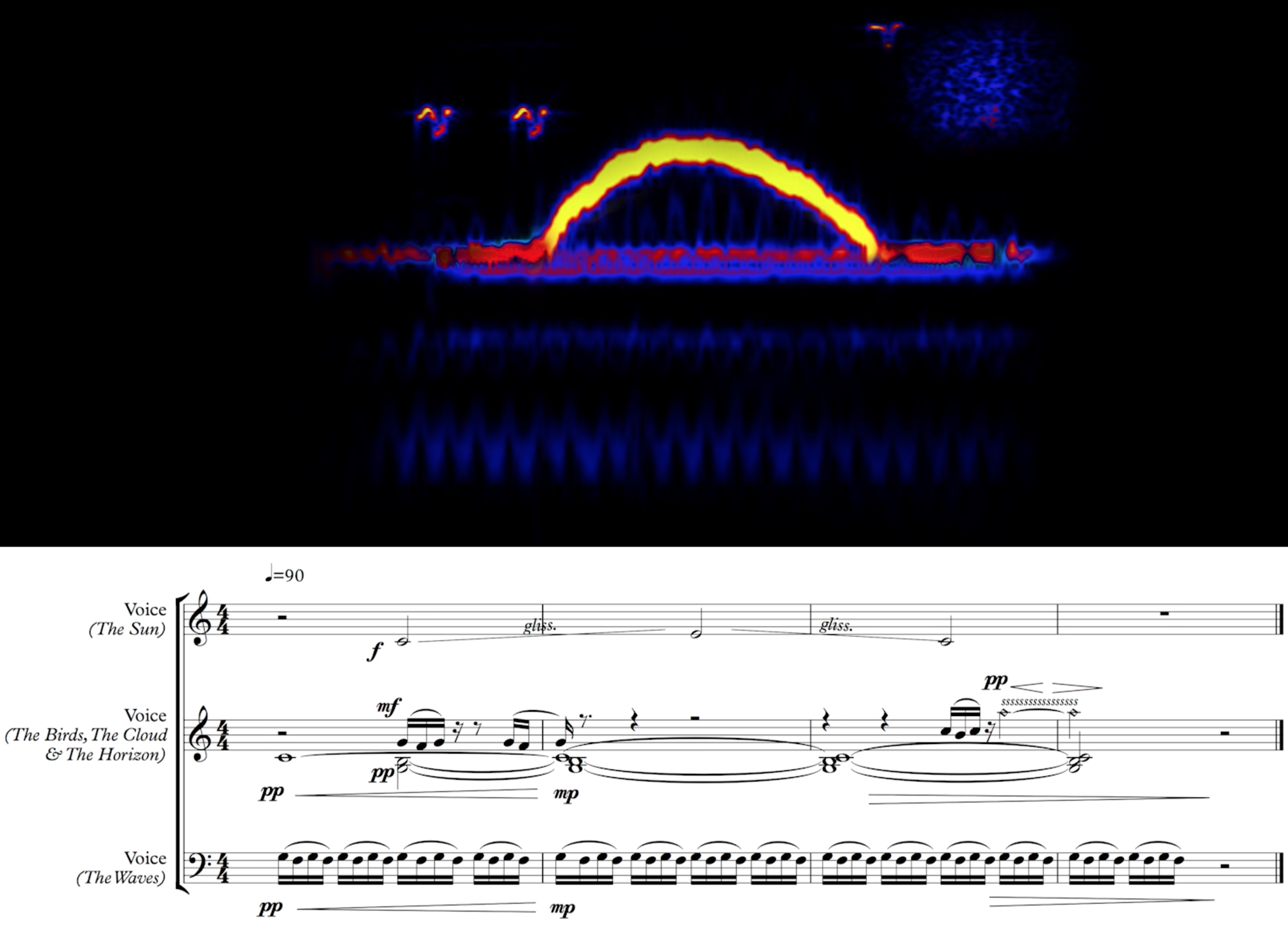}
	\caption{A drawing within a spectrogram. Image from \url{https://youtu.be/N2DQFfID6eY}.}
	\label{sunset}
\end{figure}

Translating visual cues into auditory feedback is a well-known process among digital artists. Different techniques have been previously adopted for auditory image representations.
One of the earliest approaches in digital art comes from the digital conversion of image into an audio spectrogram. 
This method involves the translation of the image's vertical components into the frequency of spectral components and the horizontal dimensions of the image in how the spectral properties of the sound evolve during time. 
An example of the application of this technique is {\em Windowlicker} by Aphex Twin (1999), to include in his piece the image-sound conversion of a spiral (see Figure \ref{Windowlicker}).
A similar approach, that instead uses traditional instruments, has been adopted in an experiment conducted by Classical Music Reimagined: the sunset is translated into music (in the form of a score), performed and recorded, and then translated back into an image (in both cases using the spectrogram method), see Figure \ref{sunset}.
	
Research between machine learning and meaningful artistic applications is a fertile field of research \citep{sturm}. 
Rather than adopting a one-to-one relationship between images and sounds, we propose a possible computational application of the method discussed in Section \ref{mathematical details} using machine learning approaches.
Specifically, we aim to utilize image recognition and segmentation process using the TensorFlow API \citep{abadi}. 
TensorFlow is an open source machine learning framework, originally developed by Google to empower a wide range of users with easy deployment of computation in different platforms.

\begin{figure}
\centering
\includegraphics[width=6cm]{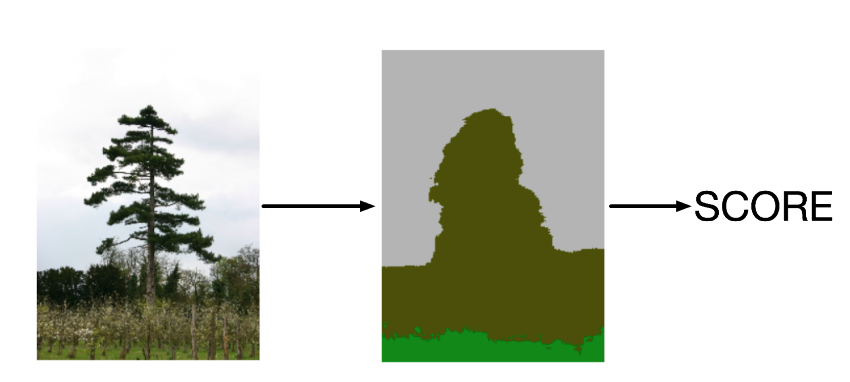}
\caption{A schematic procedure: the input is an image; then, the envelope if extracted, and it is translated into a musical score. A more precise score takes into account also patterns of the image. The pictures are adapted from \cite{vladlen}.}
\label{tree}
\end{figure}

\begin{figure}
\centering
\includegraphics[width=8cm]{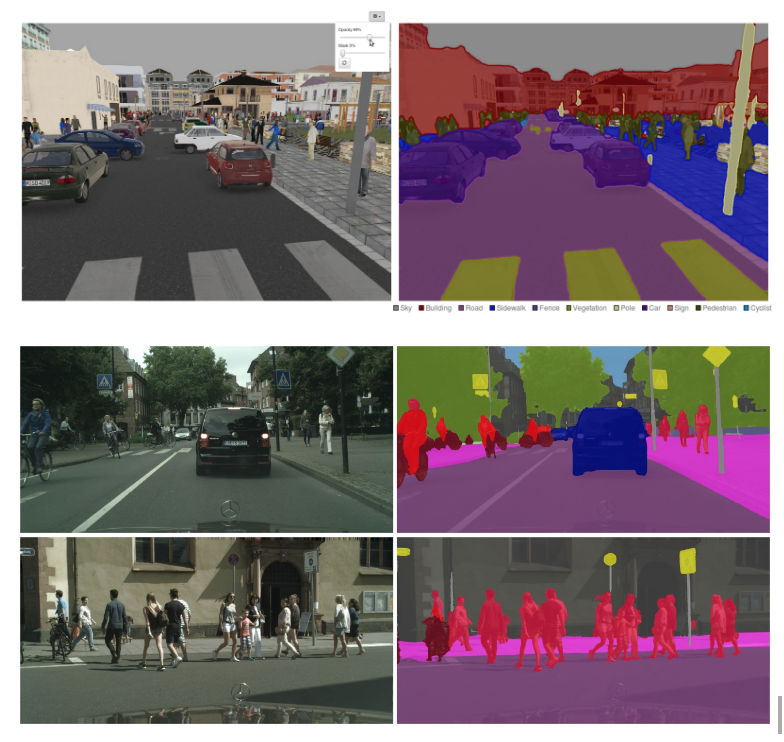}
\caption{An example of envelope extraction from different images as input.
Images from \cite{heinrich}, top, and \cite{vladlen2}, bottom.}
\label{traffic}
\end{figure}

By using TensorFlow's image recognition functionality, we will be able to identify the main content of the image, i.e., a tree in the case of Figure \ref{tree}, and a street with cars in the case of Figure \ref{traffic}. TensorFlow can be used to apply machine learning to both music and images.
The image recognition exploits a segmentation algorithm implemented using TensorFlow, that segments the image.
At this point, the application can draw contours and mutate these into musical parts.
The degree of segmentation may deal with more precise or coarse-graining segmentation of the same object, but also with higher and lower level objects: for example a detailed view of a mountain will require trees, a detailed view of trees will require branches and leaves, and so on. Finally, leaves have veins that are similar in shape to branches, with a fractal-like self-similarity. Mathematically, we may describe this via nested categories \citep{maclane, knots}. In fact, this reminds us of a n-categorical depiction. Categories can be used in programming, suggesting new approaches. This topic may strengthen the connections between mathematical, artistic, and philosophical views of the reality around us, within an overall aesthetic perspective.

\subsection{Details of the pseudocode}
This research informs the implementation of a software that empowers musicians in creating the skeleton of their compositions based on images.
Precisely, it can generate a MIDI file out of an image, which can then be converted into a sound file or a musical score.
We report the pseudocode of this future software; see Algorithm \ref{alg1}.
As a general remark, the main differences between our pseudocode and implementations such as Xenakis' UPIC tablet \citep{xenakis} lie in the discussed role of time, in the separation between patterns and envelopes, in the potential applicability of the method to tridimensional images, and in the use of the conjecture of gestural similarity.

As input, the algorithm takes (i) an image to convert, (ii) the pitch range within which the music has to be realized, (iii) the chosen instrumentation, and (iv) the discretization level. This last one corresponds to the resolution of the image analysis and the number of musical features (pitch progression, interval, dynamic, and pauses) that have to be considered in the image-music conversion process.

After passing the check on the image file format (pseudocode line \ref{alg:imageCheck}), the image recognition process takes place (line \ref{alg:recognitionProcess}), otherwise the user is asked to select a different image. The image recognition process uses a deep convolutional neural network, implemented using TensorFlow, to distinguish different objects within the image.
If objects are detected, macro and micro aspects are extracted. Macro aspects refer to the objects contour (line \ref{alg:contourAnalysis}).
Micro aspects are extracted by decomposing the image contour into smaller objects (see Section \ref{mathematical details}), which are then segmented into patterns made of dots, straight lines, and curved lines (line \ref{alg:segmentation}).

Subsequently, the pattern-music relationships are defined with a series of tests.
The first three tests establish the technique with which the music is to be played and how to define pitch and time characteristics.
Referring to concepts described through Figure \ref{basis}, lines \ref{alg:legato}, \ref{alg:softStaccato}, and \ref{alg:violentStaccato} of the pseudocode test if the music is to be played with {\em legato}, soft {\em staccato}, or violent {\em staccato} techniques respectively.
If a pattern is made of curved lines, then the notes will be indicated as {\em legato}, if the pattern is made of dots, as soft {\em staccato}, and if made of straight lines, as violent {\em staccato}. Precisely, the longer the lines, the more violent is the {\em staccato}.
Afterwards, pitch progressions (line \ref{alg:progression}), interval (\ref{alg:interval}), dynamics (\ref{alg:dynamic}) and pauses between notes (\ref{alg:pauses}) are established.
Pitch progression has a direct relationship with (visual) pattern directions (upward/downward), which is useful to establish if pitches are rising or lowering. The interval between consecutive notes is defined by looking at the pattern angle: the higher the angle, the bigger is the interval between notes. Dynamics are determined by pattern thickness -- the thicker the line in the image to be sonified, the louder the notes. Finally, pauses are defined by analyzing white spaces between each element of the pattern -- the wider the white spaces, the longer the pauses.
Notes created by analyzing patterns are then musically organized to represent the macro aspects of the given image, that is, the whole image contour, its envelope (line \ref{alg:objectEnvelopePatterns}).
As the last step, the musical score is transposed to the chosen musical instrument and register, and the MIDI file is generated.

Our pseudocode is contextualized within machine learning studies. It is worth  mentioning machine-learning studies are more and more being re-read in light of categories, using compositionality and functor definitions. This is true of Bayesian machine learning \citep{culbertson} and of backpropagation seen as a functor \citep{fong}. However, to the best of our knowledge, our approach potentially joining music, image, machine learning, and categories together is something novel and it can enter into dialogue with other ongoing or already established studies.
\begin{algorithm}[h]
\begin{footnotesize}
	\caption{MIDI file generator}\label{alg1}
	\begin{algorithmic}[1]
		\State \textbf{Input:} total time, pitch range, discretization level, instrument selection
		\State \textbf{Output:} musical score
		\If{valid image is entered} \label{alg:imageCheck}
			\State object recognition \label{alg:recognitionProcess} 
			\If{objects are recognized}
				\State objects contour analysis related to discretization level \label{alg:contourAnalysis}
				\For{each object}
					\State object segmentation $\rightarrow$ patterns generation related to discretization level  \label{alg:segmentation}
						\If{pattern is a curved line} \label{alg:legato}
							\State musical technique $\rightarrow$ legato
						\EndIf
						\If{pattern is composed of dots} \label{alg:softStaccato}
							\State musical technique $\rightarrow$ soft staccato
						\EndIf
						\If{pattern is a straight line} \label{alg:violentStaccato}
							\State musical technique $\rightarrow$ violent staccato
						\EndIf
						\If{pitch progression feature is selected} \label{alg:progression}
							\State pattern direction (up/down) $\rightarrow$ pitch progression (low/high)
						\EndIf
						\If{interval feature is selected} \label{alg:interval}
							\State pattern angle (small/big) $\rightarrow$ interval (small/big)
						\EndIf
						\If{dynamic feature is selected} \label{alg:dynamic}
								\State pattern thickness (thin/thick) $\rightarrow$ sound dynamic (piano/forte)
						\EndIf
						\If{blank space between patterns} \label{alg:pauses}
							\State space wideness (narrow/wide) $\rightarrow$ pause length (short/long)
						\EndIf
					\State notes transposition along object contour \label{alg:objectEnvelopePatterns}
				\EndFor
				\State notes transposition according to selected musical instrument \label{alg:transposition}
				\State create MIDI file \label{alg:output}
			\EndIf
			\State error: no objects have been recognized. Please, upload a different image
		\EndIf
		\State error: the uploaded image cannot be processed. Please, upload a different image
	\end{algorithmic}
	\end{footnotesize}
\end{algorithm}



\section{Discussion: open questions and a possible cognitive experiment}\label{experiment}


The proposed approach can be seen as a method to compose music and to analyze music and images. The aesthetic output can raise new questions about the strength of connections between music and images, and the depth of our understanding of these relationships.
In particular, future research would involve specific cases of multi-modal perception, choice of the bases, and experimental validation of the minimal-shape principle cited in Section \ref{mathematical details}.
In particular, two main questions arise:
\begin{enumerate}
\item How to find a good basis for images, general enough to be adapted in different cases, and where each of the basis elements can be meaningfully mapped into sound via gestural similarity? 
\item With this method of basis change with well-chosen visual kets and visual/sound basis elements' correspondences, is it possible/easier to translate a visual Gestalt into a musical Gestalt?
\end{enumerate}



In a future experiment, the participants can be given a few different complex forms and several possibilities for simple shapes (geometric curves, points, lines) as visual kets, and then be asked to select the minimum number of simple shapes to approximate each complex form, making the superposition recognizable as a simplified version of the original complex form.
For each complex form, we can evaluate how many simple shapes are needed, and which ones they are. This can help establish a principle of {\em minimum} and to recognize which shapes are more needed than others.
Then, each simple shape should be associated with a simple sound sequences. The second step of the experiment should test the appropriateness of this association. Participants might be asked to associate sonified complete forms with each visual forms, rating their similarity. In fact, the outcomes of a future experiment can also suggest which restrictive conditions may be needed to get an effective sonification and a Gestalt translation. In designing such a new experiment, we could refer to a recent experimental study, where participants were asked to draw lines and make arm gestures according to musical sequences; examples involved  vocal sequences belonging to different genres and styles \citep{kelkar}.  


Further conceptual developments of the proposed work can deal with the dialogue between continuous and discrete, and with new applications of Deleuze's studies within a Gestalt framework.
The duality between discretization and continuity is implied in Deleuze's studies about time, not only in music. Deleuze made a distinction between ``striated time'' and ``smooth time.'' Striated time to mean a discrete time, in opposition with a smooth time, a continuous one. We may see the first as a quantitative, measured time, while the second as a qualitative time. A similar distinction appears in Pierre Boulez's lessons in Darmstadt, with a distinction between pulsed and un-pulsed time \citep{boulez}. We can subdivide time in small parts, but some elements allow us to perceive it as a smooth time in a Gestalt perspective. A discrete sequence of events in time, as a sequence of photographic shots, versus a continuous flow, as a continuous movement, is relevant also in musical rendition of visual shapes. In Deleuze's philosophical research, musical conjectures have been applied to other areas of human knowledge, and vice versa. This strengthens the relevance of interdisciplinary studies, that can be adopted as metaphors in several other areas. Our open questions about the translatability of Gestalt can be connected with the problem of ``translation'' in the generalized sense of ``translatability of knowledge.'' This is the object of recent philosophical research \citep{alunni_translation}, that also involves diagrammatic thinking \citep{alunni}. Our work can contribute to the debate about possible shared, crossmodal elements not only in the arts, but in perception: according to \citealt[page 62]{gestalt_music}, this also deals with ``the question of the existence of universals of knowing and perception.''
Also, this can remind one of the universals in philosophy.

Finally, some new, brave questions may arise, in particular regarding not only the Gestalt-translatability, but also the translatability of {\em beauty} and of {\em beautiful Gestalt}: can a beautiful visual shape generate a beautiful music and vice versa? Should we carefully choose all parameters, including discretization levels, to control resulting harmonies, and so on, to translate beauty? We expect that the answers is yes, we have to carefully play around with the degrees of freedom allowed by the proposed method, which primarily aims only to preserve the recognizability of the original shape. This could open a debate about hypothetic quantitative definitions of beauty independently from a specific artistic language. We wonder if the perception of {\em beauty} follows some criterion of equilibrium between {\em anticipation} and {\em surprise}, that can be quantitatively defined in terms of pattern variety levels --- with some {\em aesthetic indeterminacy principle}. The relationship between beautiful Gestalt and gestural similarity should also be investigated. The connection between gestural similarity and expressivity in music has already been discussed \citep{mannone}. The concept of beauty can be lifted above a specific artistic language; it is the case of smoothness in visual lines and in singing, a visual analogy for a good vocal technique used by Truslit \citep{repp}. Finally, the interest of further studies about Gestalt invariance can be useful not only within an artistic framework, but also for practical, and specifically diagnostic purposes: it is the case of the sonification of  patients' brain, to create, in the mind of the listener, an auditory idea of such a complex 3D-physical object \citep{roginska}.


\section{Conclusions}\label{conclusions}

{In this article, we started with insights from philosophy, cognition, physics, and mathematical theory of music, looking for correlations among these different areas. Mathematics and diagrammatic thinking with their abstraction power can be used to formalize music and mapping strategies from image to sound. We expected to be able to decompose complex images and musical sequences in terms of simple elements combined to be together with suitable kets and coefficients. Following the Deleuzian {\em rhizomatic thought}, it is possible to look for links among these elements of images and sounds. 

With the idea of a quantization of images and sounds, we borrowed the formalism of Dirac notation from quantum physics, and we connected each step through diagrams. We formally described a method to derive music from visual shapes based on the {\em gestural similarity} criterion between simple shapes and simple musical sequences. This method can be used as a compositional technique to develop creativity and sonify different shapes and forms. However, this research also raised open-ended questions about the translatability of Gestalt from a sensory domain to another. Such a question requires cognitive experiments to help improve our theoretical model. We also described a pseudocode that could automatize (almost) each passage of the described method. Further developments of the proposed research may involve computer science applications, psychological validation of the method, and evaluation of the validity criteria for a hypothetical Gestalt translation. 
These criteria would give new blood to mathematical research in this framework with the definition of conservation principles.

The ultimate goal, in a STEAM (Science, Technology, Engineering, Arts, and Mathematics) perspective, would be a better connection and a dream of unification of different areas of human knowledge.

\section*{Authors' contributions}
Maria Mannone developed the initial idea, the mathematical section, and contributed to discussion and conclusions. Federico Favali contributed to the conceptual introduction. Balandino Di Donato contributed to the pseudocode section. Luca Turchet contributed to the cognitive application and helped revise the initial manuscript. All the authors exchanged ideas and gave feedback on each other's work.

\section*{Acknowledgements}


The authors are grateful to the physicist Peter beim Graben and to the mathematicians Giuseppe Metere and Olivia Caramello for comments on the mathematical section.
The authors thank the physicist Lucia Rizzuto for reading the manuscript and for her insightful questions, included the role of space separation in image-music translations and the role of {\em beauty}.
The authors also thank the engineer and computer scientist Davide Rocchesso for his comments and suggestions about crossmodal perception.
The authors thank the anonymous reviewers for their comments, as well as the mathematicians Emmanuel Amiot and Jason Yust for their careful work as editors.




\end{document}